# H-Harmonic Maaß-Jacobi Forms of Degree 1
## The Analytic Theory of Some Indefinite Theta Series

Martin Westerholt-Raum[1]

It was shown in previous work that the one-variable $\widehat{\mu}$-function defined by Zwegers (and Zagier) and his indefinite theta series attached to lattices of signature $(r+1,1)$ are both Heisenberg harmonic Maaß-Jacobi forms. We extend the concept of Heisenberg harmonicity to Maaß-Jacobi forms of arbitrary many elliptic variables, and produce indefinite theta series of "product type" for non-degenerate lattices of signature $(r+s,s)$. We thus obtain a clean generalization of $\widehat{\mu}$ to these negative definite lattices. From restrictions to torsion points of Heisenberg harmonic Maaß-Jacobi forms, we obtain harmonic weak Maaß forms of higher depth in the sense of Zagier and Zwegers. In particular, we explain the modular completion of some, so-called degenerate indefinite theta series in the context of higher depth mixed mock modular forms. The structure theory for Heisenberg harmonic Maaß-Jacobi forms developed in this paper also explains a curious splitting of Zwegers's two-variable $\widehat{\mu}$-function into the sum of a meromorphic Jacobi form and a one-variable Maaß-Jacobi form.

∎ real-analytic Jacobi forms
∎ generalized $\widehat{\mu}$-functions
∎ mixed mock modular forms

MSC Primary   11F50
MSC Secondary 11F27

In his celebrated thesis [Zwe02], Zwegers employed the so-called $\widehat{\mu}$-function to provide an automorphic completion of the until then mysterious mock theta functions. The $\widehat{\mu}$-function is a real-analytic Jacobi form of one modular and two elliptic variables. A remarkable fact was commented on by Zagier in [Zag09]: The "two-variable" $\widehat{\mu}$-function can be written as the sum of a meromorphic Jacobi form and a real-analytic Jacobi form that only depends on the difference of the two elliptic variables[2].

$$\widehat{\mu}(\tau,u,v) = \frac{\zeta(\tau,u) - \zeta(\tau,v) - \zeta(\tau,u-v)}{\theta(\tau,u-v)} + \widehat{\mu}(\tau,u-v)$$

where $\zeta$ is the Weierstrass $\zeta$-function, and $\theta$ is the Jacobi $\theta$-function. The second term is the "one-variable" $\widehat{\mu}$-function, which we denote, abusing notation, by the same letter as the original $\widehat{\mu}$-function. One outcome of the present paper is a natural explanation for this behavior of $\widehat{\mu}$. The construction of $\widehat{\mu}$ can be naturally phrased in terms of indefinite theta series. We extend this construction to more general lattices.

Classical Jacobi forms were defined in [EZ85], and have been applied in many contexts since then. In some cases, the generating functions of interesting arithmetic quantities turn out to be Jacobi forms [GZ98; Zag91]; in other cases, classical Jacobi forms and their generalizations have been used to understand the structure of modular forms of different type. For example, Jacobi forms occur as Fourier-Jacobi coefficients of holomorphic and non-holomorphic Siegel modular forms—see [BRR12b; Koh94] for an explanation of how Fourier-Jacobi coefficients can be obtained from the latter. Jacobi forms also serve as a tool to better understand elliptic modular forms, quasimodular forms [Zag94], and mock modular forms. Quasimodular forms, for example, occur as Taylor coefficients of classical Jacobi forms. A more recent accomplishment that is based on Jacobi forms, and which is closely connected to the subject of this paper, are the findings by Zwegers in [Zwe02]. He provided three different ways, all based on non-classical Jacobi forms, to understand mock modular forms—see [Ram00] for details on the latter. First, he defined the $\widehat{\mu}$-function, a real-analytic Jacobi form which specializes at certain torsion points to automorphic completions of mock theta functions.

---

[2]Zagier's formula has a sign mistake, as one referee pointed out. Following that referee's suggestion, we sketch a proof in Example 2.10

Second, he defined indefinite theta series for lattices of signature $(r-1,1)$, which are also real-analytic Jacobi forms. They can be employed in a similar way as the $\widehat{\mu}$-function to understand mock theta functions. Third, Zwegers analyzed Fourier coefficients of meromorphic Jacobi forms, in order to obtain mock modular forms.

Motivated by this success of real-analytic Jacobi forms (defined in an ad-hoc way), several attempts were made to give a precise definition of real-analytic Jacobi forms and, more specifically, harmonic weak Maaß-Jacobi forms. In the past few years, several such definitions, all based on the Casimir operator for the extended real Jacobi group, were suggested by Berndt and Schmidt, Pitale, Bringmann and Richter, Conley and the author, and Bringmann, Richter and the author [BR10; BRR12b; BS98; CR10; Pit09]. In order to discuss these definitions, recall that Jacobi forms are functions $\phi : \mathbb{H} \times \mathbb{C}^l \to \mathbb{C}$, depending on a modular variable $\tau \in \mathbb{H} \subset \mathbb{C}$ in the Poincaré upper half plane and elliptic variables $z \in \mathbb{C}^l$. The index of a Jacobi form is an $l \times l$ matrix. A Jacobi form is semi-holomorphic if it is holomorphic as a function of $z$. The Casimir operator is a certain invariant, central differential operator that annihilates constant functions.

Berndt and Schmidt, and Pitale gave definitions of real-analytic Jacobi forms that were motivated by representation theoretic ideas, therefore restricting themselves to functions that satisfy a polynomial growth condition with respect to the modular variable. By their definition, a real-analytic Jacobi form is an eigenfunction of the Casimir operator. In addition, Berndt and Schmidt require a real-analytic Jacobi form to be an eigenfunction of another differential operator which is invariant, but *not* central, and which is similar to $\Delta^{\mathrm{H}}$ defined in [BRR12b]. This is elaborated on in more detail in Section 1. Pitale then showed that it suffices to consider semi-holomorphic forms in order to study smooth vectors in autormorphic representations for the extended Jacobi group. This led him to require that Maaß-Jacobi forms be semi-holomorphic eigenfunctions of the Casimir operator.

The work by Bringmann and Richter introduced a new idea. Restricting to functions that are annihilated by the Casimir operator, they relaxed the growth condition, requiring at most exponential growth, and called the Jacobi forms that arise this way harmonic Maaß-Jacobi forms. In order to distinguish them from the real analytic Jacobi forms mentioned so far, we will call them *harmonic weak Maaß-Jacobi forms*. It is important to notice that functions that only satisfy a weak growth condition currently cannot be incorporated into a satisfactory representation theoretic framework. However, they gain importance by the tremendous amount of applications in which harmonic weak Maaß forms [BF04] and harmonic weak Maaß-Jacobi forms show up – see, for example, [Bru02; DIT11; DMZ12]. Note that, even though Bringmann and Richter formally did not impose any further condition in a formal way, their work treats only the semi-holomorphic case.

The weak growth condition was used in [CR10] to give a definition of semi-holomorphic harmonic weak Maaß-Jacobi forms of lattice index. In the context of [CR10], this type of Jacobi forms is relevant because of its connection with Siegel modular forms of higher genus [Rau12]. In Section 3 of [CR10], it was shown that the vanishing conditions with respect to analogs of the above $\Delta^{\mathrm{H}}$, in this setting, lead to semi-holomorphic functions, if the Jacobi index is not scalar. This is an important observation, which is treated in more detail in Proposition 2.6 of the present paper.

The later work [BRR12b] focused on scalar Jacobi indices. When allowing certain kinds of singularities, the class of harmonic weak Maaß-Jacobi forms that are annihilated by $\Delta^{\mathrm{H}}$ is strictly larger than the class of semi-holomorphic harmonic weak Maaß-Jacobi forms. A complete structure theory of *Heisenberg-harmonic* (H-harmonic) Maaß-Jacobi forms with scalar Jacobi indices was built up. We restrict ourselves to reminding the reader that the one-variable $\widehat{\mu}$-function is an H-harmonic Maaß-Jacobi form in the very sense of [BRR12a].

**The definition of** H**-harmonic Maaß-Jacobi forms.**   As discussed above, results on real-analytic Jacobi forms, so far, either restrict to the semi-holomorphic case or to the case of scalar Jacobi indices. The two-variable $\widehat{\mu}$-function is neither semi-holomorphic, nor has it scalar Jacobi index. In order to study it as a real-analytic Jacobi form, we develop the theory of H-harmonic Maaß-Jacobi forms of arbitrary indices $L$, that are lattices. To explain the definition, we recall in more detail the main result of Section 3 in [CR10]. Let $\mathrm{L}^{\mathrm{JH}}$ and $\mathrm{R}^{\mathrm{JH}}$ denote the lowering and raising operators with respect to the elliptic variables—explicit expressions are given in Section 1. These are maps that assign to every $b \in L \otimes \mathbb{R}$ a differential operator, and we write $\mathrm{L}^{\mathrm{JH}[b]}$ and $\mathrm{R}^{\mathrm{JH}[b]}$ for the images under $b$. We call them the lowering and raising operator in direction of $b$. In [CR10], it is shown that, if the rank of $L$ exceeds 1 and $L$ is non-degenerate, then any smooth function from the Jacobi upper half space $\mathbb{H} \times (L \otimes \mathbb{C})$ to $\mathbb{C}$ that is annihilated by $\mathrm{R}^{\mathrm{JH}[b]} \mathrm{L}^{\mathrm{JH}[b]}$ for all $b \in L \otimes \mathbb{R}$ is holomorphic in $z$. In other words, a coordinate independent definition of H-harmonic Maaß-Jacobi forms leads to semi-holomorphic forms whenever the Jacobi index is not scalar. From the perspective taken in [CR10], coordinate independence is a reasonable assumption, which holds automatically for Maaß-Jacobi forms that are obtained from real-analytic Siegel modular forms and orthogonal modular forms. The two-variable $\widehat{\mu}$-function, however, is not semi-holomorphic. One is thus led to consider coordinate dependent H-harmonicity. Fixing a basis $B$ of linear independent vectors in $L \otimes \mathbb{R}$, we consider *Heisenberg harmonic functions with* H-*harmonicity B* that, by definition, are annihilated by $\mathrm{R}^{\mathrm{JH}[b]} \mathrm{L}^{\mathrm{JH}[b]}$ for all $b \in B$—see Definition 2.1. The basic theory of H-harmonic Maaß-Jacobi forms is developed in Section 2.

**Indefinite theta series.**   It is natural to ask for examples of H-harmonic Maaß-Jacobi forms beyond Zwegers's $\widehat{\mu}$-function. We construct theta series for some indefinite lattices, following closely the definition proposed in [Zwe02] in the case of signature $(l_+, 1)$.

Given a lattice $L$ with bilinear form $\langle \cdot, \cdot \rangle_L$ of signature $(l_+, 1)$, Zwegers's indefinite theta functions depend on a pair $(c_1, c_2)$ of non-positive vectors in $L \otimes \mathbb{R}$. They provide modular completions of

$$\sum_{v \in L} \big(\mathrm{sgn}\langle c_1, v\rangle_L - \mathrm{sgn}\langle c_2, v\rangle_L\big) \exp\big(L[v]\tau + \langle v, b\rangle_L\big) = \widetilde{\sum_{\substack{v \in L \\ \mathrm{sgn}\langle c_1, v\rangle_L = -\mathrm{sgn}\langle c_2, v\rangle_L}}} \exp\big(L[v]\tau + \langle v, b\rangle_L\big),$$

where for technical reasons we introduce some suitable $b \in L \otimes \mathbb{R}$ such that the sum converges. The tilde decorating the sum refers to signs and boundary terms that we suppress to present a clearer picture of the construction. The condition $\langle c_1, v\rangle_L$ and $\langle c_2, v\rangle_L$ on the right hand side can be interpreted as restricting summation to a cone in $L \otimes \mathbb{R}$.

For arbitrary lattices there are several possibilities to restrict summation to achieve convergence. Utilizing the intersection of cones that are defined by two vectors each is the most straightforward one. Given a set $C$ of pairs of vectors $(c_1, c_2)$ that span mutually orthogonal spaces of $L \otimes \mathbb{R}$, we set

$$\begin{aligned}\mathrm{cone}_L(C) &= \big\{v \in L : \forall (c_1, c_2) \in C : \mathrm{sgn}\langle c_1, v\rangle_L = -\mathrm{sgn}\langle c_2, v\rangle_L\big\} \\ &= \bigcap_{(c_1, c_2) \in C} \big\{v \in L : \mathrm{sgn}\langle c_1, v\rangle_L = -\mathrm{sgn}\langle c_2, v\rangle_L\big\}.\end{aligned}$$

If $\#C = l_-$, where $(l_+, l_-)$ is the signature of $L$ and if all $c_1, c_2$ are negative, then

$$\widetilde{\theta}^C_L(\tau, 0) = \widetilde{\sum_{v \in \mathrm{cone}_L(C)}} \exp\big(L[v]\tau + \langle v, b\rangle_L\big) \tag{0.1}$$

converges. In Section 3, we provide their modular completion. As a special case, if we choose negative vectors $c_1$ and isotropic vectors $c_2$, we obtain H-harmonic Maaß-Jacobi forms for $B = \{c_1 : (c_1, c_2) \in C\}$.

Their image under the Heisenberg $\xi$-operator, discussed in Proposition 2.4, is a skew theta series. It is attached to the majorant of $L$ that is defined by taking the negative of span $B$. In particular, there are plenty of preimages of the same skew theta series, which are distinguished by different choices of bases $B$ of span $B$.

In case that for all $c_1, c_2 \in L \otimes \mathbb{Q}$ the intersection of $L \otimes \mathbb{Q}$ with span$\{c_1, c_2\}$ is two-dimensional, the theta series (0.1) can be written as a sum of products of theta series for lattices of signature $(1,1)$ and one for a positive definite lattice. For this reason, theta series of the kind that we treat are occasionally called theta series of product type. Note that under these specific assumptions on $C$ modular completions of (0.1) can be furnished by employing Zwegers's indefinite theta series. We will elaborate on this in Section 3.1.

From a purely philosophical standpoint, one expects that indefinite theta series of product type should show up rather frequently—many counting problems feature degenerate cases. We give two examples, which one of the referees suggested to mention. In their work on Torus knots Hikami and Lovejoy encountered indefinite theta series for which they could not provide a modular completion [HL15]. A concrete expression can be found in Theorem 5.6 of loc. cit. where the last sum runs over *three* variables each restricted with respect to its sign. That is, the cone that appears in a corresponding indefinite theta series has three walls. If this theta series is degenerate, it should be possible to express that cone as a union or suitable intersection of cones defined by two pairs of walls as in (0.1). Specifically, we expect that Hikami's and Lovejoy's theta series is, after multiplication with a definite theta series, of the form

$$\sum_{v \in \text{cone}} \exp(L[v]\tau), \quad \text{where} \quad \text{cone} = \{v \in L : \text{sgn}\langle c_1, v\rangle_L = \text{sgn}\langle c_2, v\rangle_L = \text{sgn}\langle c_3, v\rangle_L\}$$

with suitable $c_1$, $c_2$, $c_3$. The question is whether cone $= \bigcup_C \text{cone}_L(C)$ for a suitable collection $C$'s as above.

Another example of possibly degenerate theta series can be found in work of Lau and Zhou on Open Gromov-Witten Potentials [LZ14]. Formulas (4.10) and (4.11) both feature sums over cones with three walls.

**Theta-like decompositions.** Another foundation to our understanding of H-harmonic Maaß-Jacobi forms and indefinite theta series are theta-like decompositions as in [BRR12a]. The theta-like decomposition introduced for Maaß-Jacobi forms of scalar Jacobi index in [BRR12a] provides a more flexible way to construct examples. We prove Theorem 4.2, which extends the theta-like decomposition studied so far to the case of arbitrary Jacobi indices. Recall the statement in the case of scalar indices. An H-harmonic Maaß-Jacobi form $\phi$ of index $-m < 0$ can be written as

$$\sum_{l \,(\text{mod}\, 2m)} h_l(\tau)\, \widehat{\mu}_{m,l}(\tau, z) + \psi(\tau, z),$$

where the $h_l$ are the components of a vector-valued elliptic modular form, the $\widehat{\mu}_{m,l}$ are functions depending only on $m$ and $l$, and $\psi$ is a meromorphic Jacobi form. In our setting such a decomposition result must incorporate additional meromorphic terms. In that matter, it is interesting to note that Bringmann, Creutzig, Rolen, and Zwegers recently showed that the meromorphic term $\psi$ also admits a decomposition, namely by means of partial theta functions, if the Jacobi index $m$ is negative [BCR14; BRZ15].

A prototypical decomposition of a Jacobi form of index $\left(\begin{smallmatrix} 1 & 0 \\ 0 & -m \end{smallmatrix}\right)$, which can still be phrased in terms of Zwegers's $\widehat{\mu}$-function, is

$$\sum_{l \,(\mathrm{mod}\, 2m)} \psi_{m,l}(\tau, z_1)\, \widehat{\mu}_{m,l}(\tau, z_2) + \psi(\tau, z),$$

where the $\psi_{m,l}$ are meromorphic Jacobi forms, that depends on $\phi$.

The $\widehat{\mu}$-function by Zwegers has been one of the most prominent players in the field of real-analytic modular forms. It also appears in the theta-like decomposition, since its image under the Heisenberg $\xi$-operator, it was shown in [BRR12a], is a unary theta series. We generalize $\widehat{\mu}$ by employing previously constructed indefinite theta series. For a negative definite lattice $L$ and an orthogonal basis $B$ of $L \otimes \mathbb{R}$, the images of our functions $\widehat{\mu}_L^B$ under the Heisenberg $\xi$-operator are anti-holomorphic theta series attached to $L$.

**Splittings of** H-**harmonic Maaß-Jacobi forms.** Zwegers's $\widehat{\mu}$-function falls under Definition 2.1, and its splitting can be explained in the setting of H-harmonic Maaß-Jacobi forms. In Section 2.3 we prove that if $L$ is degenerate, then any H-harmonic Maaß-Jacobi form must be a sum of meromorphic and anti-meromorphic functions on $L_0 \otimes \mathbb{R}$, the totally isotropic part of $L \otimes \mathbb{R}$. For example, in the case of Jacobi index $L = \begin{pmatrix} 0 & 0 \\ 0 & -m \end{pmatrix}$ Proposition 2.9 and Theorem 4.2 imply that an H-harmonic Maaß-Jacobi form $\phi$ can be written as

$$\phi(\tau, (z_1, z_2)) = \psi(\tau, (z_1, z_2)) + \psi_{\overline{\mathrm{hol}}}(\tau, z_1)\phi_{\overline{\mathrm{hol}}}(\tau, z_2) + \psi_{\mathrm{hol}}(\tau, z_1)\phi_{\mathrm{hol}}(\tau, z_2)$$

for a meromorphic function $\psi$, anti-meromorphic and H-harmonic functions $\psi_{\overline{\mathrm{hol}}}$ and $\phi_{\overline{\mathrm{hol}}}$, and meromorphic and H-harmonic functions $\psi_{\mathrm{hol}}$ and $\phi_{\mathrm{hol}}$, respectively.

**Restrictions to torsion points.** In the spirit of, for example, [Zwe02], it is interesting to study restrictions of H-harmonic Maaß-Jacobi forms to torsion points. In fact, most of the contemporary theory on mock theta functions is formulated in terms of such restrictions. In Section 4, we describe their analytic properties. In particular, we connect H-harmonic Maaß-Jacobi forms with Zagier's and Zwegers's notion of harmonic weak Maaß forms of higher depth, or equivalently mixed mock modular forms of higher depth—see Section 4.2 for a definition[3]. Let $\mathbb{M}_k^{[d]}$ denote the space of depth $d$ harmonic weak Maaß forms. Then, for example, any H-harmonic Maaß-Jacobi form $\phi$ that is not singular at $z = 0$ gives

$$\phi(\tau, 0) \in \mathbb{M}_k^{[d]}.$$

By restricting H-harmonic Maaß-Jacobi forms to torsion points one obtains sums of products of harmonic weak Maaß forms. Such products cannot be characterized by differential operators, a paucity which let emerge the approach of mixed mock modular forms. Our results reconcile in parts the approach taken by geometers, who tend to focus on harmonic modular forms, and physicists, who often encounter products of harmonic weak Maaß forms and holomorphic modular forms as generating series. For example, characters of Kac-Moody Lie superalgebras [KW01; KW15] are typically not mock theta functions, but mixed mock modular forms, depending on the signature of the Lie superalgebra's Cartan matrix.

We suggest to study which of the mixed mock modular forms that have been encountered so far can be obtained as "holomorphic parts" of restrictions of H-harmonic Maaß-Jacobi forms to torsion points. Specifically, if there is, say, a counting problem with coefficients $c(n)$ such that

$$\sum_n c(n) \exp(2\pi i n \tau)$$

---

[3]There is no definition of higher depth harmonic weak Maaß forms in the literature, but it has been communicated by by Zagier and Zwegers in some talks

is a mixed mock modular form of higher depth, then we suggest to try to refine it in a natural way $c(n) = \sum_r c(n, r)$ such that

$$\sum_{n,r} c(n, r) \exp\big(2\pi i \, (n\tau + r(z))\big)$$

is a "mock Jacobi form". This has already been carried out in [BO15] for Kac-Wakimoto characters that were studied in [KW01]. An approach in the same spirit has also been helpful to investigate moments of partition counting functions [BMR14], where Taylor expansions of mock Jacobi forms occurred naturally.

**Acknoledgment.** *The author thanks Kathrin Bringmann, Olav Richter, and Sander Zwegers for helpful discussions and for their remarks. He is especially grateful to one of the referees who helped by his comments working around some stylistic glitches.*

# 1 Preliminaries

**§1.1 Lattices.** A lattice is a free $\mathbb{Z}$-module $L$ together with an $\mathbb{R}$-valued quadratic form $L[\,\cdot\,]$ on $L$. The rank of $L$ will be denoted by $l$. We say that $L$ is integral if $L[\,\cdot\,]$ takes values in $\mathbb{Z}$. The associated rational, real, and complex spaces are denoted by $L_\mathbb{Q} = L \otimes \mathbb{Q}$, $L_\mathbb{R} = L \otimes \mathbb{R}$, and $L_\mathbb{C} = L \otimes \mathbb{C}$, respectively. There is a bilinear form $\langle z, z'\rangle_L = L[z + z'] - L[z] - L[z']$ canonically attached to $L$. It extends to a linear form on $L_\mathbb{C}$ (i.e. a form which is complex linear in both the first and second component), which we also denote by $\langle \cdot, \cdot\rangle_L$.

A vector $v \in L$ is called isotropic if $L[v] = 0$. The maximal totally isotropic subspace $\{v \in L : \forall v' \in L : \langle v, v'\rangle_L = 0\}$ will be denoted by $L_0$, while its dimension is denoted by $l_0$. We call $L$ non-degenerate if $l_0 = 0$. Writing $(l_+, l_-)$ for the signature of $L$, where $l_+$ and $l_-$ are the dimensions of maximal positive and negative definite subspaces, we therefore have $l = l_+ + l_0 + l_-$. The abbreviation $l_\pm = l_+ + l_-$ will appear frequently. It is standard to call $L$ positive or negative semi-definite if $l_- = 0$ or $l_+ = 0$, respectively. We say that $L$ is positive or negative definite, if $l_+ = l$ or $l_- = l$. A totally isotropic lattice satisfies by definition $l_0 = l$.

Fixing an ordered basis for $L$, we can identify $L$ with a Gram matrix $2m = 2m_L \in \mathrm{Mat}_l^\mathrm{T}(\mathbb{R})$. The determinant of $2m$ is independent of any choice. We let the reduced covolume $|L|$ be the determinant of the matrix $2m_{L/L_0}$ that is associated with the non-degenerate lattice $L/L_0$.

The real dual $L_\mathbb{R}^\vee$ is defined for arbitrary $L$ and consists of all linear functions on $L_\mathbb{R}$. We call $L^\vee = \{v^\vee \in L_\mathbb{R}^\vee : \forall v \in L : v^\vee(v) \in \mathbb{Z}\}$ the dual of $L$. It can be identified with $\{v \in L \otimes \mathbb{Q} : \forall v' \in L : \langle v, v'\rangle \in \mathbb{Z}\}$, if $L$ is non-degenerate. In this case, we write $\mathrm{disc}(L)$ for the discriminant module $L^\vee/L$. If $L$ is non-degenerate, then $L^\vee \subseteq L_\mathbb{R}^\vee$ via $v^\vee(v) = \langle v^\vee, v\rangle_L \in \mathbb{Z}$ for $v^\vee \in L^\vee$ and $v \in L$. We define a scalar product $\langle \cdot, \cdot\rangle_L$ on $L_\mathbb{R}^\vee$ as follows. On $(L_0)_\mathbb{R}^\vee$ it is zero. A complement to $(L_0)_\mathbb{R}^\vee$ is given by the inclusion $(L/L_0)_\mathbb{R} \hookrightarrow L_\mathbb{R}^\vee$, on which we set $L\big[\langle v, \cdot\rangle_L\big] = 4L[v]$ for $v \in L_\mathbb{R}$. This quadratic form is the same that arises in the case of $l_0 = 0$ from the inclusion $L^\vee \supseteq L$.

The Weil representation associated to a non-degenerate lattice $L$ is a representation of the metaplectic cover $\mathrm{Mp}_2(\mathbb{Z})$ of $\mathrm{SL}_2(\mathbb{Z})$ on the group algebra $\mathbb{C}[\mathrm{disc}(L)]$ (see [Sko08]). A natural basis for $\mathbb{C}[\mathrm{disc}(L)]$ is given in terms of $\mathfrak{e}_v$, where $v$ runs through $\mathrm{disc}(L)$. By abuse of notation we write $S$ and $T$ for the generators of the metaplectic cover of $\mathrm{SL}_2(\mathbb{Z})$ that project to the corresponding generators $S = \left(\begin{smallmatrix} 0 & -1 \\ 1 & 0 \end{smallmatrix}\right)$ and $T = \left(\begin{smallmatrix} 1 & 1 \\ 0 & 1 \end{smallmatrix}\right)$ of $\mathrm{SL}_2(\mathbb{Z})$. In terms of our basis of $\mathbb{C}[\mathrm{disc}(L)]$ and these generators of $\mathrm{Mp}_2(\mathbb{Z})$, the Weil representation is defined as:

$$\rho_m(T)\mathfrak{e}_v := e(L[v])\mathfrak{e}_v, \qquad \rho_m(S)\mathfrak{e}_v := \frac{1}{\sigma_L \sqrt{2|L|}} \sum_{v' \in \mathrm{disc}(L)} e\big(-\langle v, v'\rangle_L\big)\mathfrak{e}_{v'}, \tag{1.1}$$

where $\sigma_L = \sqrt{2|L|}^{-1} \sum_{v \in \mathrm{disc}(L)} e(-L[x])$ is the signum of $L$. Here, we use the shorthand notation $e(x) = \exp(2\pi i x)$ for $x \in \mathbb{C}$ that will appear later again. Throughout the paper we will pass from $\mathrm{SL}_2(\mathbb{Z})$ to its metaplectic cover whenever necessary and without further mentioning it.

Given any set of vectors $v_1, \ldots, v_n$ we denote their span by $\mathrm{span}\{v_1, \ldots, v_n\}$. It will be clear by the context, whether we mean the span over $\mathbb{Z}$, $\mathbb{R}$, or $\mathbb{C}$. The orthogonal complement of $v_1, \ldots, v_n$ will be denoted by $\{v_1, \ldots, v_n\}^\perp$, or in the case of $n = 1$ by $v_1^\perp$. For a subset $B$ of $L_\mathbb{Q}$, the span of $B$ is denoted by $L_B = L \cap \mathrm{span}_\mathbb{Q} B$.

Given a subset $B$ of $L \otimes \mathbb{R}$, we define

$$B_+ = \{b \in B : L[b] > 0\}, \quad B_- = \{b \in B : L[b] < 0\}, \quad B_0 = \{b \in B : L[b] = 0\}. \tag{1.2}$$

**§1.2 Jacobi forms.** The Poincaré upper half plane and the Jacobi upper half space attached to a lattice $L$ are

$$\mathbb{H} = \{\tau = x + iy : y > 0\} \subset \mathbb{C}, \quad \text{and} \quad \mathbb{H}^{(JL)} = \mathbb{H} \times (L \otimes \mathbb{C}^2).$$

The latter is isomorphiphic to $\mathbb{H}^{(JI)} = \mathbb{H} \times \mathbb{C}^l$ in a non-canonical way. Typically, elements of $\mathbb{H}^{(JL)}$ are written as pairs $(\tau, z)$, where $z = u + iv$ with $u, v \in L_\mathbb{R}$. Recall the notation $e(x) = \exp(2\pi i x)$ for $x \in \mathbb{C}$. The variable $q$ stands for $e(\tau)$. Given $r \in L_\mathbb{R}^\vee$, we set $\zeta^r = e(r(z))$.

Multiplication in the real Jacobi group attached to $L$

$$G^{(JL)} := \mathrm{SL}_2(\mathbb{R}) \ltimes (L \otimes \mathbb{R}^2)$$

is given by

$$(\gamma, \lambda, \mu) \cdot (\gamma', \lambda', \mu') = \big(\gamma\gamma', (\lambda, \mu)\gamma' + (\lambda', \mu')\big).$$

Here, a typical element of $G^{(JL)}$ is denoted by $g^J = (\gamma, \lambda, \mu)$, where $\lambda, \mu \in L \otimes \mathbb{R}$. The pair $(\lambda, \mu)$ is viewed as an element of $L \otimes \mathbb{R}^2$, on which $\gamma'$ acts trivial on the first component and by its standard representation from the right on the second component. Note that $G^{(JL)}$ is independent of the quadratic form $q_L$. The discrete subgroup $\Gamma^{(JL)} = \mathrm{SL}_2(\mathbb{Z}) \ltimes (L \otimes \mathbb{Z}^2) \subset G^{(JL)}$ is called the full Jacobi group.

The action of $\mathrm{SL}_2(\mathbb{R})$ on $\mathbb{H}$ is given by

$$\gamma\tau = \frac{a\tau + b}{c\tau + d},$$

where here and throughout the element $\gamma$ of $\mathrm{SL}_2(\mathbb{R})$ is written as $\begin{pmatrix} a & b \\ c & d \end{pmatrix}$. The real Jacobi group acts on the Jacobi upper half plane via

$$(\gamma, \lambda, \mu)(\tau, z) = \Big(\gamma\tau, \frac{z + \lambda\tau + \mu}{c\tau + d}\Big).$$

Fix an integer $k$ and a lattice $L$ as above. The action of $G^{(JL)}$ on $\mathbb{H}^{(JL)}$ admits a cocyle

$$\alpha_L^J(g^J, z) = e\Big(\frac{-cL[z]}{c\tau + d} + \langle z, \lambda\rangle_L + L[\lambda]\tau\Big), \tag{1.3}$$

which leads to the following slash action on $C^\infty(\mathbb{H}^{(JL)})$.

$$\big(\phi|_{k,L}(\gamma, \lambda, \mu)\big)(\tau, z) = (c\tau + d)^{-k} \alpha_L^J(g^J, z)\, \phi\big((\gamma, \lambda, \mu)(\tau, z)\big). \tag{1.4}$$

This is the usual Jacobi slash action.

We say that a function $\phi : \mathbb{H}^{(JL)} \to \mathbb{C}$ has *non-moving singularities*, if there are finitely many linear maps $\lambda^\vee : L_\mathbb{Q} Q \to \mathbb{Q}$ and finitely many $\alpha, \beta \in \mathbb{Q}$, such that the singularities of $\phi$ are located at $\{z : \lambda(z) = \alpha + \beta\tau\} + L + \tau L$. We say that such a singularity at $(\tau_0, z_0)$ has (i) meromorphic type, (ii) almost meromorphic type, or (iii) real-analytic quotient type, if there is a neighborhood $U \subset \mathbb{H}$ of $\tau_0$, a function $z_0 : \mathbb{H} \to L \otimes \mathbb{C}$ with $z_0(\tau_0) = z_0$, and a function $\psi : U \setminus \{(\tau, z_0(\tau)) : \tau \in U\}$ which (i) is meromorphic, (ii) is the quotient of a real-analytic by a holomorphic functions, or (iii) is the quotient of two real-analytic functions such that $\phi - \psi$ can be continued to a real-analytic function in a neighborhood of $(\tau_0, z_0)$. A complex valued function on $\mathbb{H}^{(JL)}$ that is meromorphic as a function of the second component $z$ will be called semi-meromorphic. Any such function is called semi-holomorphic if it has no singularities.

*Remark 1.1.* It appears that the term "non-moving poles" was coined by Zagier. The intuition behind it is that torsion points $\alpha\tau + \beta$, $\alpha, \beta \in L_\mathbb{Q}$ should be viewed as fixed while varying $\tau$. The typical example of a modular form with non-moving singularities of meromorphic type are inverse theta series.

A meromorphic Jacobi form of weight $k$ and index $L$ with non-moving singularities is a meromorphic function $\phi : \mathbb{H}^{(JL)} \to \mathbb{C}$ such that $\phi|_{k,L} \gamma^J = \phi$ for all $\gamma^J \in \Gamma^{(JL)}$. The space of such forms is denoted by $\mathscr{M}J_{k,L}$, where $\mathscr{M}$ refers to meromorphicity.

Finally, it is clear that the concept of vector valued elliptic modular forms extends to Jacobi forms. The Weil representation, previously described as a representation of $\mathrm{SL}_2(\mathbb{Z})$ or its metaplectic double cover, can be viewed as representation of the full Jacobi group by defining $\rho_L(\lambda) = \rho_L(\mu) = \mathrm{id}$, the identity.

**§ 1.3 Differential operators.** We discuss the theory of covariant and invariant differential operators for the extended real Jacobi group. The results coincide with the ones in [CR10], if $L$ is non-degenerate. Even if $L$ is degenerate, we still employ the results of [CR10] in a crucial way. The reader is referred to this work for details on the Lie theoretic background. Note that our notation significantly differs from the one that appeared in previous work. It is close to the nowadays common notation for the classical differential operators for Maaß forms—see page 177 of[Maa64]. This notation seems thus better suited for usage outside the context of Lie groups, but inside the scope of modular forms and Jacobi forms.

Let $\mathrm{D}^\infty(\mathbb{H}^{(JL)})$ be the algebra of differential operators on $\mathbb{H}^{(JL)}$. The subalgebra of $G^{(JL)}$-invariant differential operators for the slash action defined in (1.4) is denoted by $\mathbb{D}(k,L)$. We may regard the imaginary part $v$ of $z \in L_\mathbb{C}$ as an element $\langle v, \cdot \rangle_L$ of $L_\mathbb{R}^\vee \otimes \mathrm{D}^\infty(\mathbb{H}^{(JL)})$, where $\mathrm{D}^\infty(\mathbb{H}^{(JL)})$ is the space of smooth differential operators on $\mathbb{H}^{(JL)}$. We view $\partial_z$ and $\partial_{\overline{z}}$ as elements of $L_\mathbb{R}^\vee \otimes \mathrm{D}^\infty(\mathbb{H}^{(JL)})$, so that, in particular, $\langle v, \partial_z \rangle_L \in \mathrm{D}^\infty(\mathbb{H}^{(JL)})$.

$$\mathrm{L}^{J}_{k,L} := -2iy\left(y\partial_{\overline{\tau}} + \langle v, \partial_{\overline{z}} \rangle_L\right), \quad \mathrm{R}^{J}_{k,L} := 2i\left(\partial_\tau + y^{-1}\langle v, \partial_z \rangle_L + 2\pi i\, y^{-2} L[v]\right) + k y^{-1},$$
$$\mathrm{L}^{JH}_{k,L} := -iy\partial_{\overline{z}}, \qquad\qquad \mathrm{R}^{JH}_{k,L} := i\partial_z - 2\pi y^{-1}\langle v, \cdot \rangle_L. \tag{1.5}$$

As is common by now, we suppress the subscripts $k, L$, if they are clear form the context. We write $\mathrm{L}^{JH[v]}$ and $\mathrm{R}^{JH[v]}$ for the evaluation of $\mathrm{L}^{JH}$ and $\mathrm{R}^{JH}$, respectively, at $v \in L_\mathbb{R}$.

*Remark 1.2.* The superscript J of the raising and lowering operator should remind the reader of *J*acobi forms, on which they act. The superscript H refers to the *H*eisenberg subgroup of $G^{(JL)}$. The lowering and raising operators $\mathrm{L}^{JH}$ and $\mathrm{R}^{JH}$ both act only on the elliptic variable $z$, which originates in the Heisenberg subgroup of $G^{(JL)}$.

The commutation relations of lowering and raising operators are

$$[L^J, R^J] = -k, \quad [L^J, L^{JH}] = 0, \qquad [L^J, R^{JH}] = -L^{JH},$$
$$[R^J, R^{JH}] = 0, \quad [L^{JH}, R^J] = R^{JH}, \quad [L^{JH}, R^{JH}] = -\pi \langle,\rangle_L \in L_{\mathbb{R}}^\vee \otimes L_{\mathbb{R}}^\vee. \tag{1.6}$$

They can be verified readily by means of Helgason's theory of differential operators [Hel59] as displayed in [CR10], and a computation of the images of $y$ or $v$.

The commutation relations in (1.6) show that we can view the raising and lowering operators in conjunction with one further element k as generators of an abstract algebra. The additional element k acts on Jacobi forms by multiplication with their weight. The has commutation relations of k are

$$[L^J, k] = 2L^J, \quad [R^J, k] = -2R^J, \quad [L^{JH}, k] = L^{JH}, \quad [R^{JH}, k] = -R^{JH}.$$

We write $D^J$ for this algebra, and $\mathbb{D}^J$ for the k centralizer subalgebra of $D^J$. It consists of elements that act on Jacobi forms as invariant differential operators. To emphasize dependence on $L$, we occasionally add the subscript $L$ to $D^J$ and $\mathbb{D}^J$.

We define a Casimir operator $\mathscr{C}^J_{k,L} \in \mathbb{D}^\infty(\mathbb{H}^{(JL)})$, extending the expression in [CR10] to the case of degenerate $L$. It is given by

$$-2R^J L^J + i\left(R^J \langle L^{JH}, L^{JH}\rangle_L - L^J \langle R^{JH}, R^{JH}\rangle_L\right) - \tfrac{1}{2}\left(L[R^{JH}]L[L^{JH}] - \langle R^{JH}\langle R^{JH}, L^{JH}\rangle_L, L^{JH}\rangle_L\right)$$
$$-\tfrac{i}{2}(2k-l-3)\langle R^{JH}, L^{JH}\rangle_L, \tag{1.7}$$

which equals

$$-2\Delta_{k-l_\pm/2} + \tfrac{y^2}{\pi i}\left(\partial_{\bar\tau} L[\partial_z] + \partial_\tau L[\partial_{\bar z}]\right) - 8y\partial_\tau \langle v, \partial_{\bar z}\rangle_L + \tfrac{y^2}{32\pi^2}\left(4L[\partial_{\bar z}]L[\partial_z] - \langle \partial_z, \partial_{\bar z}\rangle_L^2\right)$$
$$+ \tfrac{y}{2\pi i}\langle v, \partial_{\bar z}\rangle_L \langle \partial_z, \partial_u\rangle_L - \tfrac{(2k-l+1)y}{8\pi}\langle \partial_{\bar z}, \partial_u\rangle_L + 2\langle v\langle v, \partial_{\bar z}\rangle_L, \partial_{\bar z}\rangle_L + (2k-l_+ - l_- - 1)i\langle v, \partial_{\bar z}\rangle_L. \tag{1.8}$$

We easily check that $\mathscr{C}^J_{k,L}$ commutes with all invariant differential operators by using that this is the case for non-degenerate $L$ and that $\mathscr{C}^J$ depends continuously on $L[\cdot]$.

Suppose that $L$ is degenerate. Then differentials with respect to $z \in L_0 \otimes \mathbb{C} \subseteq L_\mathbb{C}$ do not occur in $\mathscr{C}^J$. Instead, additional operators arise from the totally isotropic part of $L$. From the formal element k, we obtain an $\mathfrak{sl}_2$-triple $(L^J, k, R^J)$. As stated before, the kernel of the commutator $[k \cdot,]$ consists of invariant differential operators. Let

$$D^{JH_0} = \mathbb{C}\left[L^{JH[v]}, R^{JH[v]} : v \in L_0 \otimes \mathbb{R}\right]$$

be the subalgebra of $D^J$ attached to the totally isotropic part of $L_\mathbb{R}$.

**Proposition 1.3.** *The center of* $D^J_L$ *is generated by* $\mathscr{C}^J$, *viewed as an element of* $D^J_L$ *by* (1.7), *and the* $\mathfrak{sl}_2$-*invariants* $H^0(\mathrm{span}\{L^J, k, R^J\}, D^{JH_0})$.

*Remark 1.4.* Note that we refer, in the above proposition, to the center of $D^J$, not the one of $\mathbb{D}^J$. The latter might be larger, but we believe that they agree.

*Proof.* We have already asserted that $\mathscr{C}^J$ is central in $\mathbb{D}^J$. By results in [CR10] it spans the center of $D^J_{L/L_0}$. Since further $D^{JH_0}$ is the kernel of $D^J_L \to D^J_{L/L_0}$ the result follows. ∎

Let us fix notation for the following extra differential operators, which lie in the part of the center of $\mathbb{D}^{\mathrm{J}}$ that arises from $L_0 \otimes \mathbb{R}$.

$$\Delta_L^{\mathrm{JH}[v,v']} := \mathrm{L}^{\mathrm{JH}[v]} \mathrm{R}^{\mathrm{JH}[v']} - \mathrm{L}^{\mathrm{JH}[v']} \mathrm{R}^{\mathrm{JH}[v]}, \qquad v, v' \in L_0 \otimes \mathbb{R}, v \neq v'. \tag{1.9}$$

The classical weight $k$ Laplace operator on $\mathbb{H}$ arises from $\mathscr{C}^{\mathrm{J}}$ when we consider the lattice $L$ of rank 0, which we provisionally denote by $\emptyset$.

$$\Delta_k := \mathrm{R}^{\mathrm{J}}_{k,\emptyset} \mathrm{L}^{\mathrm{J}}_{k,\emptyset} = 4y^2 \partial_\tau \partial_{\overline{\tau}} - 2kiy\partial_{\overline{\tau}}. \tag{1.10}$$

It factors as a product $\Delta_k = \xi_{2-k}\xi_k$, where

$$\xi_k f = y^{k-2} \overline{\mathrm{L}^{\mathrm{J}}_{k,\emptyset} f} \tag{1.11}$$

is the classical $\xi$ operator, that first appeared in [BF04].

The Heisenberg Laplace operator is similar to $\Delta_k$, but acts merely on the $z$-variable of functions on $\mathbb{H}^{(JL)}$. Given $v \in L_{\mathbb{R}}$, we set

$$\Delta_L^{\mathrm{JH}[v]} := \mathrm{R}^{\mathrm{JH}[v]}_{k-1,m} \mathrm{L}^{\mathrm{JH}[v]}_{k,m} = y \partial_z(v) \partial_{\overline{z}}(v) + 4\pi i \langle v,v \rangle_L \partial_{\overline{z}}(v). \tag{1.12}$$

**§1.4 Harmonic weak Maaß forms.** We revisit briefly the theory of harmonic weak Maaß forms. A nice exposition can be found in [Bru02]. Vector valued elliptic modular forms are invariant under the $|_{k,\rho}$ slash action, which is associated to some weight $k \in \frac{1}{2}\mathbb{Z}$ and a type, i.e. a finite dimensional, complex representation $\rho$ of $\mathrm{SL}_2(\mathbb{Z})$ or $\mathrm{Mp}_2(\mathbb{Z})$. The slash action is defined by

$$\bigl(f|_{k,\rho} \gamma\bigr)(\tau) = (c\tau+d)^{-k} \rho(\gamma)^{-1} f\bigl(\gamma\tau\bigr),$$

where $\gamma\tau$ denotes the action of $\mathrm{SL}_2(\mathbb{R})$ on $\mathbb{H}$. Recall that we pass to the metaplectic cover in case that $k \notin \mathbb{Z}$.

The space of vector-valued modular forms of weight $k$ and type $\rho$ is the space of holomorphic functions $f: \mathbb{H} \to V$ that satisfy (i) $f|_{k,\rho}\gamma = f$ for all $\gamma \in \Gamma$ and (ii) $f(\tau) = O(1)$ as $y \to \infty$. We denote this space by $\mathrm{M}_k(\rho)$. The space of weakly holomorphic elliptic modular forms is the space of holomorphic functions $f: \mathbb{H} \to V$ where the second condition is weakened: $f(\tau) = O(e^{ay})$ for some $a > 0$ as $y \to \infty$. This space is denoted by $\mathrm{M}_k^!(\rho)$.

The space of harmonic weak Maaß forms of weight $k$ and type $\rho$, which is denoted by $\mathbb{M}_k(\rho)$, consists of real-analytic functions $f: \mathbb{H} \to \mathbb{C}$ satisfying (i) $f|_{k,\rho}\gamma = f$ for all $\gamma \in \Gamma$, (ii) $\Delta_k f = 0$, and (iii) $f(\tau) = O(e^{ay})$ as $y \to \infty$ for some $a > 0$.

**§1.5 Skew Jacobi forms.** We close the preliminaries with the definition of skew(-holomorphic) Jacobi forms of matrix index. The original definition of skew Jacobi forms of scalar indices $m$ was given in [Sko90]. It was formulated referring to a slash action

$$|c\tau+d|^{-1} \overline{(c\tau+d)}^{k-1} e\Bigl(\frac{-mcz^2}{c\tau+d}\Bigr) \phi\Bigl(\frac{a\tau+b}{c\tau+d}, \frac{z}{c\tau+d}\Bigr).$$

This allowed for choosing positive Jacobi indices $m$. However, in the case of general Jacobi indices $L$, such a choice is slightly impractical, when used in conjunction with the differential operator $\mathrm{L}^{\mathrm{JH}}$, as we will see later. More importantly, it hides an interesting structural analogy that occurs when passing from real-anlytic Jacobi forms to skew Jacobi forms. For this reason we choose a different equivalent slash action to define skew Jacobi forms.

Given $L$, fix a maximal negative definite subspace $L_\mathbb{R}^{\mathrm{sk}} \subseteq L_\mathbb{R}$. Since $L_\mathbb{R}^{\mathrm{sk}}$ is non-degenerate, it allows for an orthogonal decomposition of $z$ into $z_1$ and $z_2 \in L_\mathbb{R}^{\mathrm{sk}} \otimes \mathbb{C}$. We write $\lambda = \lambda_1 + \lambda_2$, $\mu = \mu_1 + \mu_2$, and $g_1^{\mathrm{J}} = (g, \lambda_1, \mu_1)$, $g_2^{\mathrm{J}} = (g, \lambda_2, \mu_2)$. For functions $\phi : \mathbb{H}^{(\mathrm{J} L)} \to \mathbb{C}$, the skew slash action attached to $L_\mathbb{R}^{\mathrm{sk}}$ is defined by

$$\phi \big|_{k,L}^{\mathrm{sk}[L_\mathbb{R}^{\mathrm{sk}}]} g^{\mathrm{J}} = |c\tau + d|^{-l_-} (c\tau + d)^{l_-/2 - k} \alpha_L^{\mathrm{J}}(\gamma_1^{\mathrm{J}}, z_1) \overline{\alpha_L^{\mathrm{J}}(\gamma_2^{\mathrm{J}}, z_2)} \big(\phi \circ (\gamma, \lambda, \mu)\big). \tag{1.13}$$

The cocycle $\alpha^{\mathrm{J}}$ can be found in (1.3).

Set $L^{\mathrm{sk}} = L \cap L_\mathbb{R}^{\mathrm{sk}}$, which is by means of the restriction of $\langle \cdot, \cdot \rangle_L$, a negative definite lattice. The skew heat operator attached to $L$ is the complex conjugate of the usual one:

$$\mathbb{L}_L^{\mathrm{sk}} = \partial_{\overline{\tau}} + \tfrac{1}{8\pi i} L[\partial_{\overline{z}}]. \tag{1.14}$$

It annihilates anti-holomorphic theta series

$$\theta_{L^{\mathrm{sk}}, \overline{\nu}}^{\mathrm{sk}} = \sum_{\nu \in \overline{\nu} + L^{\mathrm{sk}}} \overline{q}^{L[\nu]} \overline{\zeta}^{\langle z, \nu \rangle_L} \tag{1.15}$$

that can be defined for $\overline{\nu} \in (L^{\mathrm{sk}} \otimes \mathbb{Q}) / L^{\mathrm{sk}}$:

$$\mathbb{L}_L^{\mathrm{sk}} \theta_{L^{\mathrm{sk}}, \overline{\nu}}^{\mathrm{sk}}(\tau, z) = \sum_{\nu \in \overline{\nu} + L^{\mathrm{sk}}} \big(-2\pi i L[\nu] + \tfrac{2}{8\pi i}(2\pi i)^2 L[\langle \nu, \cdot \rangle_L]\big) \overline{q}^{L[\nu]} \overline{\zeta}^{\langle z, \nu \rangle_L},$$

where we have explicitly written the inclusion of $\nu \in L_\mathbb{R}$ into $L_\mathbb{R}^\vee$. Recall that by definition we have $L[\langle \nu, \cdot \rangle_L] = 4 L[\nu]$, which implies $\mathbb{L}_L^{\mathrm{sk}} \theta_{L^{\mathrm{sk}}, \overline{\nu}}^{\mathrm{sk}} = 0$.

**Definition 1.5.** Let $L$ be a lattice, and let $L_\mathbb{R}^{\mathrm{sk}}$ and $z = z_1 + z_2$ be as above. A function $\phi : \mathbb{H}^{(\mathrm{J} L)} \to \mathbb{C}$ that is real-analytic except for non-moving singularities of real-analytic quotient type, is called a skew Jacobi form of weight $k$, index $L$, and skewness $L_\mathbb{R}^{\mathrm{sk}}$ if it satisfies the following conditions:

(i) For all $\gamma^{\mathrm{J}} \in \Gamma^{(\mathrm{J} L)}$, we have $\phi \big|_{k,L}^{\mathrm{sk}[L_\mathbb{R}^{\mathrm{sk}}]} \gamma^{\mathrm{J}} = \phi$.

(ii) The function $\phi$ is meromorphic in $z_1$ and antimeromorphic in $z_2$.

(iii) We have $\mathbb{L}_L^{\mathrm{sk}} \phi = 0$.

(iv) We have $\phi(\tau, z) = O(1)$ as $y \to \infty$ if $z = \alpha \tau + \beta$, $\alpha, \beta \in L_\mathbb{R}$ is a non-singular point of $\phi$.

We will denote the space of such skew Jacobi forms by $\mathscr{M} \mathrm{J}_{k,L}^{\mathrm{sk}[L^{\mathrm{sk}}]}$.

Singularities of skew Jacobi forms cannot be located in arbitrary position. The first instance of a related discussion appears in [BRR12a], which we extend to the case of lattice indices.

**Proposition 1.6.** *Let $\phi \in \mathscr{M} \mathrm{J}_{k,L}^{\mathrm{sk}[L^{\mathrm{sk}}]}$. If $\phi$ is singular along $\{(\tau, z) : \nu^\vee(z) = \alpha \tau + \beta\}$ for $\nu^\vee \in L_\mathbb{Q}^\vee$ and $\alpha, \beta \in \mathbb{Q}$, then $\nu^\vee$ lies in the orthogonal complement of $\langle L_\mathbb{R}^{\mathrm{sk}}, \cdot \rangle_L \subseteq L_\mathbb{R}^\vee$.*

*Proof.* We prove the proposition by contradiction. Suppose that there is $\nu^{\mathrm{sk}} \in L_\mathbb{R}^{\mathrm{sk}}$ with $\nu^\vee(\nu^{\mathrm{sk}}) \neq 0$ along which some skew Jacobi form $\phi$ has singularities.

We employ coordinates $z_1, \ldots z_l$ with respect to which $L$ has the form $L[z] = \sum_i s_i z_i^2$ with $s_i \in \{-1, 0, 1\}$. Further, we may assume that $L_\mathbb{R}^{\mathrm{sk}} = \mathrm{span}\{z_1, \ldots, z_{l^{\mathrm{sk}}}\}$ for $l^{\mathrm{sk}} = \dim L_\mathbb{R}^{\mathrm{sk}}$. Since $\phi$ is anti-meromorphic with respect to $z_i$, $1 \le i \le l^{\mathrm{sk}}$, we have a local Laurent expansion

$$\sum_{n=(n_1,\ldots,n_l)\in\mathbb{Z}^l} c_n(\tau) \prod_{i=1}^{l^{\mathrm{sk}}} \overline{(z_i - z_{i,0}(\tau))}^{n_i} \prod_{i=l^{\mathrm{sk}}}^{l} (z_i - z_{i,0}(\tau))^{n_i}$$

where $z_{i,0}$ is a linear in $\tau$, and $c_n = 0$ for small enough $n_i$. By definition, $\phi$ is annihilated by the skew heat operator, which implies that

$$\sum_n \prod_{i=1}^{l^{\mathrm{sk}}} \overline{(z_i - z_{i,0}(\tau))}^{n_i} \prod_{i=l^{\mathrm{sk}}}^{l} (z_i - z_{i,0}(\tau))^{n_i}$$
$$\left( \partial_{\overline{\tau}} c_n(\tau) + \sum_{i=1}^{l^{\mathrm{sk}}} \frac{-\partial_{\overline{\tau}} \overline{z_{i,0}(\tau)}}{\overline{z_i - z_{i,0}(\tau)}} + \sum_{\substack{i,j=1 \\ i \ne j}}^{l^{\mathrm{sk}}} \frac{n_i n_j}{\overline{(z_i - z_{i,0}(\tau))(z_j - z_{j,0}(\tau))}} + \sum_{i=1}^{l^{\mathrm{sk}}} \frac{n_i(n_i+1)}{\overline{(z_i - z_{i,0}(\tau))}^2} \right) \quad (1.16)$$

vanishes. We show in three steps that $n_i \ge 0$ if $c_n \ne 0$ and $1 \le i \le l^{\mathrm{sk}}$. Exploiting symmetry, it suffices to consider the case $i = 1$. Fix $n$ with $c_n \ne 0$ that is minimal with respect to the ordering

$$n' < n \Leftrightarrow \exists i_0 : n_{i_0} < n'_{i_0} \wedge \forall i < i_0 \, n_i = n'_i. \quad (1.17)$$

First, consider the last term in (1.16), which must vanish if $n_1 < 0$, since the exponent of $z_i - z_{i,0}(\tau)$ that arises from it is minimal with respect to (1.17). From this we see that $n_1 \ge -1$. Second, if $n_1 = -1$, we inspect the next to last term in (1.16). It shows that $n_i = 0$ for $2 \le i \le l^{\mathrm{sk}}$. Third, we inspect the second term. If $z_{1,0}$ is constant, it is zero, so that the first term forces $n_1 \ge 0$. Otherwise, it yields the lowest term of (1.16), and hence again $n_1 \ge 0$. ∎

The previous proposition is essential, since it allows us to obtain theta decompositions with respect to $L^{\mathrm{sk}} \subseteq L$ for *all* skew Jacobi forms, as opposed to holomorphic Jacobi forms, for which singularities can occur.

**Proposition 1.7.** *Fix $L$ and $L^{\mathrm{sk}}$ as above, and suppose that $L^{\mathrm{sk}} \subseteq L$ splits off as a direct summand: $L = L' \oplus L^{\mathrm{sk}}$ for some lattice $L'$. Recall the decomposition of $z$ into $z_1$ and $z_2$. For any $\phi \in {}^{\mathcal{M}}\mathrm{J}_{k,L}^{\mathrm{sk}[L^{\mathrm{sk}}]}$ we have a theta decomposition*

$$\phi(\tau, z) = \sum_{v_0 \in \mathrm{disc}(L^{\mathrm{sk}})} \psi_{v_0}(\tau, z_1) \theta_{L^{\mathrm{sk}}, v_0}^{\mathrm{sk}}(\tau, z_2), \quad (1.18)$$

*for a vector valued meromorphic Jacobi form*

$$(\psi_{v_0})_{v_0 \in \mathrm{disc}(L^{\mathrm{sk}})} \in {}^{\mathcal{M}}\mathrm{J}_{k-l_-/2, L'}(\rho_{L^{\mathrm{sk}}}).$$

*Proof.* We only provide a sketch of the proof, and leave verification of the details to the reader. Fixing a skew Jacobi form $\phi$ as in the statement, we have a local Fourier expansion at $(\tau_0, z_0) \in \mathbb{H}^{(JL)}$ that are of the form

$$\sum_{\substack{n \in \mathbb{Z} \\ r \in L^\vee}} c_{(\tau_0, z_0)}(n, r; y, v) q^n \zeta^r.$$

Proposition 1.6 says that given $z \in L_{\mathbb{C}}$ the function $z_2 \in L_{\mathbb{C}}^{\mathrm{sk}} \mapsto \phi(\tau, z + z_2)$ has no singularities. In conjunction with invariance of $\phi$ under the action of $L^{\mathrm{sk}} \otimes \mathbb{Z}^2 \subset \Gamma^{(JL)}$, we find that

$$c_{(\tau_0, z_0)}(n, r; y, v) = c_{(\tau_0, z_0)}(n', r + 2r^{\mathrm{sk}}; y, v)$$

for all $r^{\mathrm{sk}} \in L^{\mathrm{sk}}$ and with $4n + L[r] = 4n' + L[r + 2r^{\mathrm{sk}}]$. From this, one concludes directly that there is a theta decomposition of the asserted form. ∎

*Remark 1.8.* We assume that $L^{\mathrm{sk}}$ splits off from $L$, but this restriction is not essential. We can always pass to suplattices of $L$ that splits into orthogonal sums, and employ vector valued Jacobi forms. This procedure requires the extend full Jacobi group defined in [CR10].

## 2 H-**harmonic Maaß-Jacobi forms**

In [CR10], Conley and the author provided a definition of H-harmonic Maaß-Jacobi forms for non-degenerate lattice indices. We chose a coordinate independent definition, and were able to show that if the rank of $L$ exceeds 1, then all instances of H-harmonic Maaß-Jacobi forms are semi-meromorphic. In other words, for such lattices, H-harmonic Maaß-Jacobi forms in the sense of [CR10] fall under the scope of previous considerations.

We now provide a definition of H-harmonic Maaß-Jacobi forms that depends on the choice of a subset $B$ of $\mathbb{P}(L_{\mathbb{R}})$. Elements of $B$ will be frequently identified with preimages in $L_{\mathbb{R}}$ to simplify notation. This section contains the study of basic properties of H-harmonic Maaß-Jacobi forms. The definition is formulated for arbitrary $B$, but we will see in Propositions 2.6 and 2.11, and Corollary 2.7 that it suffices to consider a finite sets of mutually orthogonal, negative vectors. This explains our choice of notation: The reader should think of $B$ as a basis for a negative definite subspace of $L_{\mathbb{R}}$ (very much related to $L_{\mathbb{R}}^{\mathrm{sk}}$ in the previous section).

**Definition 2.1.** Let $L$ be a lattice, and $B \subseteq \mathbb{P}(L_{\mathbb{R}})$ with span $B = L_{\mathbb{R}}$. Let $\phi : \mathbb{H}^{(JL)} \to \mathbb{C}$ be real-analytic except for non-moving singularities of real-analytic quotient type. We say that $\phi$ is an H-harmonic Maaß-Jacobi form of weight $k$, index $L$, and H-harmonicity $B$ if it satisfies the following conditions:

(i) For all $\gamma^J \in \Gamma^{(JL)}$, we have $\phi|_{k,L} \gamma^J = \phi$.

(ii) The function $\phi$ is annihilated by the center of $\mathbb{D}_{k,L}^J$. In particular, we have $\mathscr{C}_{k,L}^J \phi = 0$.

(iii) We have $\Delta^{\mathrm{JH}[b]} \phi = 0$ for all $b \in B$.

(iv) The growth condition $\phi(\tau, \alpha\tau + \beta) = O(e^{ay})$ as $y \to \infty$ is satisfied for some $a > 0$ provided that $z = \alpha\tau + \beta$, $\alpha, \beta \in L_{\mathbb{R}}$ is a non-singular point of $\phi$.

We denote the space of all such H-harmonic Maaß-Jacobi forms by $\mathscr{MJ}_{k,L}^{\Delta,\mathrm{H}[B]}$ Notation for spaces of vector valued Jacobi forms is analogous to the one for elliptic modular forms, introduced below: We refer to a representation $\rho$ of $\Gamma^{(JL)}$ in parenthesis.

*Remark 2.2.* (1) We assume in Definition 2.1 that H-harmonic Maaß-Jacobi forms have singularities of real-analytic quotient type. However, it follows from the theory that we develop in this paper, that they automatically have almost meromorphic singularities. This is a consequence of Theorem 4.2, which in this regard is based on Proposition 1.7.

(2) The definition extends to half-integral weights and indices, and to complex representations of $\Gamma^{(JL)}$. We do not treat these cases explicitly, as they lead to additional technical difficulties and do not yield further insight. Zwegers's $\widehat{\mu}$-function, however, strictly speaking does not fall under Definition 2.1. Confer [BRR12a] for more details on the latter.

(3) Equally well, one can define H-harmonic Maaß-Jacobi forms for any $\Gamma = \Gamma' \ltimes L^2$ where $\Gamma' \subset \mathrm{SL}_2(\mathbb{Z})$ has finite index.

(4) An analogous definition can be made for the skew Jacobi slash action, which subsumes the image of $\xi^{\mathrm{JH}[b]}$ defined below.

**§2.1 Various spaces of H-harmonic Maaß-Jacobi forms.** In analogy with the ideas that were presented in [BRR12a], we define several subspaces of $\mathcal{M}\mathbb{J}_{k,L}^{\Delta,\mathrm{H}[B]}$, and remark on how they might be studied. It should be clear that not all of them appear in the present work, but it seems advantageous to suggest uniform notation for sequels to this paper.

We will very soon restrict our attention to $\mathcal{M}\mathbb{J}_{k,L}^{\mathrm{H}[B]}$, which is defined by

$$\mathcal{M}\mathbb{J}_{k,L}^{\mathrm{H}[B]} = \mathcal{M}\mathbb{J}_{k,L}^{\delta,\mathrm{H}[B]} := \mathcal{M}\mathbb{J}_{k,L}^{\Delta,\mathrm{H}[B]} \cap \ker \mathrm{L}^{\mathrm{J}}.$$

This space is very much analogous to the one primarily studied in [BRR12a]. We will deduce a theta-like decomposition for the forms in it in Section 4. Notation defined here for general H-harmonic Maaß-Jacobi forms, transfers directly to those annihilated by $\mathrm{L}^{\mathrm{J}}$.

If $L$ is non-degenerate, and $B \subseteq L_{\mathbb{R}}$ does not span $L_{\mathbb{R}}$, then we set

$$\mathcal{M}\mathbb{J}_{k,L}^{\Delta,\mathrm{H}[B]} := \mathcal{M}\mathbb{J}_{k,L}^{\Delta,\mathrm{H}[B+B^{\perp}]} \cap \ker_{b \in B^{\perp}} \mathrm{L}^{\mathrm{JH}[b]}.$$

We refrain from extending this definition to degenerate lattices, since $L_0 \otimes \mathbb{R} \subset B^{\perp}$ for every $B$, so that notation would become ambiguous. Instead, we set, for $B_{\mathrm{hol}} \subseteq \mathbb{P}(L_{\mathbb{R}})$ with $B_{\mathrm{hol}} + \mathrm{span}\, B = L_{\mathbb{R}}$

$$\mathcal{M}\mathbb{J}_{k,L}^{\Delta,\mathrm{h}[B_{\mathrm{hol}}]\mathrm{H}[B]} := \mathcal{M}\mathbb{J}_{k,L}^{\Delta,\mathrm{H}[B+B_{\mathrm{hol}}]} \cap \ker_{b \in B_{\mathrm{hol}}} \mathrm{L}^{\mathrm{JH}[b]}.$$

In accordance with notation in [BRR12a], we write $\mathcal{M}\mathbb{J}_{k,L}^{\Delta,\mathrm{h}}$ for $\mathcal{M}\mathbb{J}_{k,L}^{\Delta,\mathrm{h}[L_{\mathbb{R}}]} = \mathcal{M}\mathbb{J}_{k,L}^{\Delta\mathrm{h}[L_{\mathbb{R}}]\mathrm{H}[\varnothing]}$.

*Remark 2.3.* (1) The theory of Maaß-Jacobi forms in $\mathcal{M}\mathrm{J}_{k,L}^{\Delta,\mathrm{h}}$ was developed in [CR10]. It resembles strongly the theory of harmonic weak Maaß forms.

(2) A theory similar to the one developed in [BRR12a] could be developed for $\mathcal{M}\mathrm{J}^{\Delta,\mathrm{H}[B]}$ by considering restrictions.

(3) If $L$ is non-degenerate, then $\mathcal{M}\mathbb{J}_{k,L}^{\mathrm{h}} = \mathcal{M}\mathrm{J}_{k,L}$ is the space of meromorphic Jacobi forms with non-moving singularities of meromorphic type.

**§2.2 The Heisenberg $\xi$-operators.** As in the case of $l = 1$, dealt with in [BRR12a], there is a Heisenberg $\xi$-operator, which we define now. Given $v \in L_{\mathbb{R}}$, we set

$$\xi_L^{\mathrm{JH}[b]} := \frac{\mathrm{sgn}(L[b]\langle v,b\rangle_L)}{\sqrt{|L[b]|y}} \exp\Big(\frac{-\pi \langle v,b\rangle_L^2}{L[b]\, y}\Big) \mathrm{L}_L^{\mathrm{JH}[b]}. \tag{2.1}$$

If $B$ is a set of orthogonal vectors, non of which is isotropic, then we set $\xi^{\mathrm{JH}[B]} = \prod_{b \in B} \xi^{\mathrm{JH}[b]}$.

**Proposition 2.4.** *Fix a subspace $L_{\mathbb{R}}^{\mathrm{sk}} \subseteq L_{\mathbb{R}}$ that contains no totally isotropic vector. Given a set $B \subset \mathbb{P}(L_{\mathbb{R}})$ of mutually orthogonal vectors that are also orthgonal to $L_{\mathbb{R}}^{\mathrm{sk}}$, and non of which is isotropic, we let $L_{\mathbb{R}}^{\mathrm{sk}} + B$ be the span of $L_{\mathbb{R}}^{\mathrm{sk}}$ and $B$. Then for any smooth function $\phi$ on $\mathbb{H}^{(\mathrm{J}\,L)}$ and any $g^{\mathrm{J}} \in G^{(\mathrm{J}\,L)}$, we have*

$$\xi_L^{\mathrm{JH}[B]}\bigl(\phi|_{k,L}^{\mathrm{sk}[L_{\mathbb{R}}^{\mathrm{sk}}]}\, g^{\mathrm{J}}\bigr) = \bigl(\xi_L^{\mathrm{JH}[B]}\phi\bigr)|_{k,L}^{\mathrm{sk}[L_{\mathbb{R}}^{\mathrm{sk}}+B]}\, g^{\mathrm{J}}.$$

*Proof.* It suffices to treat the case $B = \{b\}$. Write $z_b = b\langle z,b\rangle_L / 2L[b]$ and $\lambda_b = b\langle \lambda, b\rangle_L / 2L[b]$. Using known covariance of $L^{\mathrm{JH}}$ and $y^{\frac{-1}{2}}$, we only need to establish that

$$\exp\Bigl(\frac{-\pi \langle v,b\rangle_L^2}{L[b]\, y}\Bigr) e\Bigl(\frac{-cL[z_b + \lambda_b \tau]}{c\tau + d} + \langle z_b, b\rangle_L + L[\lambda_b]\tau\Bigr)$$
$$= e\Bigl(\frac{-cL[\overline{z}_b + \lambda_b \overline{\tau}]}{c\overline{\tau} + d} + \langle \overline{z}_b, b\rangle_L + L[\lambda_b]\overline{\tau}\Bigr)\Bigl(\exp\Bigl(\frac{-\pi \langle v,b\rangle_L^2}{L[b]\, y}\Bigr)\Bigr)_{\tau \mapsto \frac{-1}{\tau},\, z \mapsto \frac{z}{\tau}}.$$

Checking this for $g^{\mathrm{J}} = \bigl(\bigl(\begin{smallmatrix} 0 & -1 \\ 1 & 0 \end{smallmatrix}\bigr), 0, 0\bigr)$ and for $g^{\mathrm{J}} = \bigl(\bigl(\begin{smallmatrix} 1 & 0 \\ 0 & 1 \end{smallmatrix}\bigr), \lambda, 0\bigr)$ separately establishes the statement. ∎

*Remark 2.5.* The image of $\xi^{\mathrm{JH}[B']}$ with $B' \subseteq B$ applied to $\mathcal{M}\mathbb{J}_{k,L}^{\Delta, \mathrm{H}[B]}$, in general, is an instance of skew H-harmonic Maaß-Jacobi forms, referred to in Remark 2.2(4). If $B' = B$, we obtain proper skew Jacobi forms as in Definition 1.5. This is a consequence of Proposition 2.15.

**§2.3 Non-trivial H-harmonicities are orthogonal.** Definition 2.1 of H-harmonic Maaß-Jacobi forms refers to a set $B \subseteq \mathbb{P}(L_{\mathbb{R}})$. However, it suffices to study sets of mutually orthogonal vectors. The argument that we employ to show this is similar to the one in [CR10]. But the degenerate case does not allow for a quite straightforward generalization.

**Proposition 2.6.** *For any $B \subseteq \mathbb{P}(L_{\mathbb{R}})$, we have*

$$\mathcal{M}\mathbb{J}_{k,L}^{\mathrm{H}[B]} \subseteq \mathcal{M}\mathbb{J}_{k,L}^{\mathrm{h}[B_{\neq 0}]\mathrm{H}[B_{=0}]}, \quad \text{where} \quad \begin{aligned} B_{=0} &= \bigl\{b' \in B : \forall b \in B, b \neq b' : \langle b, b'\rangle_L = 0\bigr\}, \\ B_{\neq 0} &= \bigl\{b' \in B : \exists b \in B, b \neq b' : \langle b, b'\rangle_L \neq 0\bigr\}. \end{aligned} \qquad (2.2)$$

*Proof.* Fix $b, b' \in B$ with $b \neq b'$ and $\langle b, b'\rangle_L \neq 0$. We lift both to elements of $L_{\mathbb{R}}$, denoting these again by $b$ and $b'$. In case that $b - b' \in L_0 \otimes \mathbb{R}$, then replace, without loss of restriction, $b'$ by $2b'$. In particular, we may assume that $b - b'$ is not totally isotropic.

We have to show that any $\phi \in C^\infty(\mathbb{H}^{(\mathrm{J}\,L)})$ that is annihilated by $\Delta^{\mathrm{JH}[b]}$ and $\Delta^{\mathrm{JH}[b']}$ is already annihilated $L^{\mathrm{JH}[b]}$ and $L^{\mathrm{JH}[b']}$. Now, the remainder of the proof is the same as the proof of Theorem 3.4 of [CR10]. For convenience of the reader, we reproduce a variation of it. It suffices to consider functions $a(y,v)$ that are annihilated by $\Delta^{\mathrm{JH}[b]}$ and $\Delta^{\mathrm{JH}[b']}$. The commutator

$$\bigl[\Delta^{\mathrm{JH}[b]}, \Delta^{\mathrm{JH}[b']}\bigr] = \pi \langle b, b'\rangle_L \bigl(\mathrm{R}^{\mathrm{JH}[b']}\mathrm{L}^{\mathrm{JH}[b]} - \mathrm{R}^{\mathrm{JH}[b]}\mathrm{L}^{\mathrm{JH}[b']}\bigr)$$

also annihilates $a(y,v)$, which implies that

$$\partial_v(b')\, a(y,v) = \frac{\langle v, b'\rangle_L}{\langle v, b\rangle_L}\, \partial_v(b)\, a(y,v).$$

Because

$$\partial_v(b') \frac{\langle v, b'\rangle_L}{\langle v, b\rangle_L} = \frac{\langle b', b'\rangle_L \langle v, b\rangle_L - \langle b, b'\rangle_L \langle v, b'\rangle_L}{\langle v, b\rangle_L^2},$$

we find that $\Delta^{\mathrm{JH}[b']}$ acts on $a(y,v)$ as

$$\frac{\langle v, b'\rangle_L}{\langle v, b\rangle_L} \partial_v(b)^2 + \left(\frac{-2\pi}{y} \frac{\langle v, b'\rangle_L^2}{\langle v, b\rangle_L} + \frac{\langle b', b'\rangle_L \langle v, b\rangle_L - \langle b, b'\rangle_L \langle v, b'\rangle_L}{\langle v, b\rangle_L^2}\right) \partial_v(b).$$

We now consider the action of $\Delta^{\mathrm{JH}[b']} - \frac{\langle v, b'\rangle_L}{\langle v, b\rangle_L} \Delta^{\mathrm{JH}[b]}$, which is

$$\left(\frac{-2\pi}{y} \frac{\langle v, b'\rangle_L^2}{\langle v, b\rangle_L} + \frac{\langle b', b'\rangle_L \langle v, b\rangle_L - \langle b, b'\rangle_L \langle v, b'\rangle_L}{\langle v, b\rangle_L^2} + \frac{2\pi}{y} \frac{\langle v, b'\rangle_L}{\langle v, b\rangle_L} \frac{\langle v, b'\rangle_L^2}{\langle v, b\rangle_L}\right) \partial_v(b).$$

Therefore $a(y, v)$ is constant in direction of $b$ if the first factor is not zero. Since this factor is a nontrivial polynomial in $\langle v, b\rangle_L$ and $\langle v, b'\rangle_L$, this implies the statement. ∎

**Corollary 2.7.** *Let $B \subseteq \mathbb{P}(L \otimes \mathbb{R})$ be a set of mutually orthogonal vectors that spans a maximal nondegenerate subspace of $L_\mathbb{R}$. Then*

$$\mathcal{M}\mathbb{J}_{k,L}^{\mathrm{H}[B]} = \mathcal{M}\mathbb{J}_{k,L}^{\mathrm{h}[B_0]\mathrm{H}[B_\pm]},$$

*where $B_0 = \{b \in B : L[b] = 0\}$ and $B_\pm = \{b \in B : L[b] \ne 0\}$ are the subsets of isotropic and non-isotropic vectors in $B$.*

*Proof.* By the assumption that span $B$ is maximal non-degenerate, for each $b \in B_0$ there is $b' \in B$ with $\langle b, b'\rangle_L \ne 0$. Proposition 2.6 therefore implies the corollary. ∎

The totally isotropic subspace of $L_\mathbb{R}$ would pose particular technical problems, if only vanishing with respect to the Casimir operator is imposed. However, Definition 2.1 also features the degenerate central invariant differential operators $\Delta^{\mathrm{JH}[v,v']}$ given in (1.9). They enforce reasonable analytic behavior on $\mathbb{H} \times (L_0 \otimes \mathbb{C}) \subseteq \mathbb{H}^{(JL)}$.

**Lemma 2.8.** *Any orthogonal set $B \subseteq \mathbb{P}(L_\mathbb{R})$ with span $B = L_\mathbb{R}$ contains a generating set of $L_0 \otimes \mathbb{R}$.*

*Proof.* Passing to the quotient $L/L_0$, observe that span $B = (L/L_0) \otimes \mathbb{R}$. Further, vectors in $L/L_0 \otimes \mathbb{R}$ are orthogonal to each other if and only if they are so in $L_\mathbb{R}$. Thus we see that $B$ contains exactly $l_+ + l_-$ vectors that are not totally isotropic. The remaining vectors must span $L_0 \otimes \mathbb{R}$, since $B$ spans $L_\mathbb{R}$. ∎

**Proposition 2.9.** *Any $\phi \in \mathcal{M}\mathbb{J}_{k,L}^{\mathrm{H}[B]}$ is either holomorphic or anti-holomorphic with respect to $L_0 \otimes \mathbb{C} \subseteq L_\mathbb{C}$.*

*Proof.* Lemma 2.8 tells us that $B_0$ generates $L_0 \otimes \mathbb{R}$. By assumption, $\phi$ is annihilated by $\Delta^{\mathrm{JH}[b]}$ for $b \in B_0$. If $\phi$ is constant with respect to every $b \in B_0$, then we are done. Assuming this is not the case, we can fix one $b$ so that either $\mathrm{L}^{\mathrm{JH}[b]}\phi \ne 0$ or $\mathrm{R}^{\mathrm{JH}[b]}\phi \ne 0$. Further, fix $v \in L_0 \otimes \mathbb{R}$ which is not a multiple of $b$.

Since both $b$ and $v$ are totally isotropic, $\mathrm{L}^{\mathrm{JH}[v]} = -iy\partial_{\overline{z}}(v)$, $\mathrm{R}^{\mathrm{JH}[v]} = i\partial_z(v)$, $\mathrm{L}^{\mathrm{JH}[b]} = -iy\partial_{\overline{z}}(b)$, and $\mathrm{R}^{\mathrm{JH}[b]} = i\partial_z(b)$ commute. Therefore, $\phi$ vanishing under $\Delta^{\mathrm{JH}[v]}$ and $\Delta^{\mathrm{JH}[b]}$ is equivalent to $\phi$ being locally a sum of a holomorphic and an anti-holomorphic function in direction of $v$ and $b$, each.

Consider the equation $\Delta^{\mathrm{JH}[b,v]}\phi = 0$ which is implied by Condition (ii) of Definition 2.1 in conjunction with (1.9). It means that

$$\mathrm{L}^{\mathrm{JH}[v]} \mathrm{R}^{\mathrm{JH}[b]} \phi = \mathrm{R}^{\mathrm{JH}[v]} \mathrm{L}^{\mathrm{JH}[b]} \phi.$$

The left hand side is holomorphic with respect to $b$, the right hand side is anti-holomorphic with respect to $b$. Since $\phi$ is a sum of a holomorphic and an anti-holomorphic function in $b$, this implies that both the left and the right hand side are zero. Therefore, if $\mathrm{L}^{\mathrm{JH}[v]}\phi \ne 0$, then $\mathrm{R}^{\mathrm{JH}[b]}\phi = 0$, and if $\mathrm{R}^{\mathrm{JH}[v]}\phi \ne 0$, then $\mathrm{L}^{\mathrm{JH}[b]}\phi = 0$. This proves the claim. ∎

For reference, we fix notation for spaces that are meromorphic and anti-meromorphic on $L_0 \otimes \mathbb{R}$

$$\begin{aligned}
\mathcal{M}\mathbb{J}^{\mathrm{h}[B_{\mathrm{hol}}]\,\mathrm{h}_0\,\mathrm{H}[B]} &= \mathcal{M}\mathbb{J}^{\mathrm{h}[B_{\mathrm{hol}}]\mathrm{H}[B+L_0\otimes\mathbb{R}]} \cap \ker_{b\in L_0\otimes\mathbb{R}} \mathrm{L}^{\mathrm{JH}[b]}, \\
\mathcal{M}\mathbb{J}^{\mathrm{h}[B_{\mathrm{hol}}]\,\overline{\mathrm{h}}_0\,\mathrm{H}[B]} &= \mathcal{M}\mathbb{J}^{\mathrm{h}[B_{\mathrm{hol}}]\mathrm{H}[B+L_0\otimes\mathbb{R}]} \cap \ker_{b\in L_0\otimes\mathbb{R}} \mathrm{R}^{\mathrm{JH}[b]}.
\end{aligned} \quad (2.3)$$

*Example 2.10.* Recall that $\zeta(\tau,z)$ is the Weierstass $\zeta$ function, and $\theta(\tau,z)$ is the Jacobi theta function.

$$\theta(\tau,z) = \sum_{n\in \frac{1}{2}+\mathbb{Z}} (-1)^{n-\frac{1}{2}} q^{n^2} \zeta^n = q^{\frac{1}{8}} \zeta^{\frac{1}{2}} \prod_{n=1}^{\infty} (1-q^n)(1-\zeta q^n)(1-\zeta^{-1}q^{n-1}) \quad \text{and}$$

$$\zeta(\tau,z) = \frac{\partial_z \theta(\tau,z)}{\theta(\tau,z)} = \frac{1}{z} + \sum_{\substack{z_0 \in \mathbb{Z}+\tau\mathbb{Z} \\ z_0 \neq 0}} \left( \frac{1}{z-z_0} + \frac{1}{z_0} + \frac{z}{z_0^2} \right).$$

From this representation of $\theta$ ad $\zeta$ we directly read off their zeros and poles.

The splitting of Zwegers's two-variable $\widehat{\mu}$-function is well-know, but in [Zag09] Zagier only remarked on it in a footnote. Zwegers's definition of $\widehat{\mu}$ is

$$\widehat{\mu}(\tau,z_1,z_2) = \frac{e^{\pi i z_1}}{\theta(\tau,z_2)} \sum_{n\in\mathbb{Z}} \frac{(-1)^n e\left(\frac{n^2+n}{2}\tau + nz_2\right)}{1-e(n\tau+z_1)}$$
$$- \frac{i}{2\sqrt{\pi}} \sum_{n\in \frac{1}{2}+\mathbb{Z}} \left( -\mathrm{sgn}(n) - H^{\mathrm{H}[(-1,1)]}(y, u_1-u_2; n) \right)(-1)^{n+\frac{1}{2}} e\left(\frac{-n^2}{2}\tau + n(z_1-z_2)\right), \quad (2.4)$$

which we obtain after unreveling Formula (1.8) in [Zwe02]. Its index, Theorem 1.11 of [Zwe02] states, is the matrix $m = \frac{1}{2}\begin{pmatrix} -1 & 1 \\ 1 & -1 \end{pmatrix}$. Note that here we fix a basis for the corresponding lattice with Gram matrix $m$.

The one-variable $\widehat{\mu}$-function was never defined by Zwegers in a formal way, and Zagier in [Zag09] simply asserted existence of a one-variable $\widehat{\mu}$-function. So we have to prove that

$$\widehat{\mu}(\tau,z_1,z_2) - \frac{\zeta(\tau,z_1) - \zeta(\tau,z_2) + \zeta(\tau,z_1-z_2)}{\theta(\tau,z_1-z_2)} \quad (2.5)$$

transforms like a Jacobi form and only depends on $z_1-z_2$.

The two variable $\widehat{\mu}$ function has simple poles at $z_2 \in \mathbb{Z}+\tau\mathbb{Z}$ (from the prefactor of the first term of (2.4)) and at $z_1 \in \mathbb{Z}+\tau\mathbb{Z}$ (from the denominator inside the sum of the first term of (2.4)). The residue at $z_1 = 0$ is easily seen to be $\theta(\tau,z_2)^{-1}$, and the symmetry of $\widehat{\mu}(\tau,z_1,z_2)$ stated in Theorem 1.11 (3) of [Zwe02] shows that the residue at $z_2 = 0$ is $\theta(\tau,z_1)^{-1}$. The poles of $\zeta(\tau,z)$ lie at $z \in \mathbb{Z}+\tau\mathbb{Z}$. One thus quickly checks that the function in (2.5) has poles at $z_1-z_2 \in \mathbb{Z}+\tau\mathbb{Z}$.

We conveniently inspect the modular properties of $\zeta(\tau,z_1)-\zeta(\tau,z_2)+\zeta(\tau,z_1-z_2)$ by means of [Rol15]— also confer [Obe12]—, which implies it is a Jacobi form of index 0. The quotient by $\theta(\tau,z_1-z_2)$ yields a meromorphic Jacobi form of index $m$, the same as the one of $\widehat{\mu}(\tau,z_1,z_2)$. Therefore, (2.5) displays a real analytic Jacobi form of index $m$ whose poles are supported on $z_1-z_2 \in \mathbb{Z}+\tau\mathbb{Z}$.

The image of $\widehat{\mu}(\tau,z_1,z_2)$ under $\xi^{\mathrm{JH}[b_2]}$ vanishes, from which we infer that $\widehat{\mu}$ is meromorphic in direction of $b_2 = (1,1)$. And so is (2.5). From our above considerations of poles and residues, we further find that it is holomorphic is direction of $b_2$. Fixing $z_1-z_2 \in \mathbb{C} \setminus (\mathbb{Z}+\tau\mathbb{Z})$, and letting vary $z_1+z_2$, we obtain a holomorhpic Jacobi form of index $m[(1,1)] = 0$. By the classical theory of Jacobi forms—see, for example, Theorem 1.2 of [EZ85]—it is constant. This establishes the splitting of the two-variable $\widehat{\mu}$-function.

**§2.4 Fourier expansions.** Fourier expansions are crucial to understand further properties of H-harmonic Maaß-Jacobi forms. Even more so, since they allow us to reason about restrictions. Any term of the Fourier expansion is indexed by $n \in \mathbb{R}$ and $r \in L_\mathbb{R}^\vee$. The overall situation is similar to the case of scalar Jacobi indices: The space of Fourier coefficients $c(n, r; y, v) q^n \zeta^r$ that are annihilated by $\mathscr{C}^{\mathrm{J}}$ and $\Delta^{\mathrm{JH}[b]}$ for $b \in B$, $B \subseteq \mathbb{P}(L_\mathbb{R})$ is finite dimensional, if $\mathrm{span}\, B = L_\mathbb{R}$.

Maaß-Jacobi forms with non-moving singularities admit a local Fourier expansion as is explained in [BRR12a]. For $\phi \in \mathscr{M}\mathbb{J}_{k,L}^{\mathrm{H}[B]}$, we denote it by

$$\phi(\tau, z) = \sum_{n \in \mathbb{Z}, r \in L^\vee} c(\phi; n, r; y, v)\, q^n \zeta^r. \tag{2.6}$$

Note that, if $L$ is degenerate, then $L^\vee$ cannot be canonically identified with a subspace of $L_\mathbb{Q}$, as is common otherwise.

Propositions 2.6 and 2.9 allow us to focus on $B \subseteq \mathbb{P}(L_\mathbb{R})$ that consists of mutually orthogonal, negative vectors, and $B_{\mathrm{hol}} \subseteq \mathbb{P}(L_\mathbb{R})$ which is orthogonal to $B$ and contains no isotropic vectors. Our goal is to describe (2.6) explicitly for

$$\phi \in \mathscr{M}\mathbb{J}_{k,L}^{\mathrm{h}[B_{\mathrm{hol}}]\,\mathrm{h}_0\,\mathrm{H}[B]} \quad \text{and for} \quad \phi \in \mathscr{M}\mathbb{J}_{k,L}^{\mathrm{h}[B_{\mathrm{hol}}]\,\overline{\mathrm{h}}_0\,\mathrm{H}[B]}.$$

We split $z \in L_\mathbb{C}$ into $z_{\mathrm{hol}} \in \mathrm{span}\, B_{\mathrm{hol}}$, $z_0 \in L_0 \otimes \mathbb{C}$, and $z_b = \langle z, b\rangle_L b / 2L[b] \in \mathbb{C}b$ for $b \in B$ such that $z = z_{\mathrm{hol}} + z_0 + \sum_b z_b$.

Besides the holomorphic terms, two building blocks suffice to explicitly describe Fourier expansions. The function $H$, up to slight modifications, appeared in [BF04] first. It is usually employed to describe Fourier expansions of harmonic weak Maaß forms. The function $H^{\mathrm{H}[b]}$, adopted from [BRR12a], is adjusted to the setting of lattices indices. It features the lower incomplete gamma function $\gamma(s, y) = \int_0^y t^{s-1} e^{-t}\, dt$. We set

$$H(y; D) := e^{-y} \int_{-2y}^\infty e^{-t} t^{-k+l_\pm/2}\, dt\Big|_{y=\pi Dy/2|L|}, \quad \text{if } D \neq 0, |L| \neq 0; \tag{2.7}$$

$$H(y; D) := y^{-k+l_\pm/2}, \quad \text{if } D = 0; \text{ and} \tag{2.8}$$

$$H^{\mathrm{H}[b]}(y, v; r) := \mathrm{sgn}\big(\langle v + ry, b\rangle_L\big)\, \gamma\Big(\tfrac{1}{2}, -\pi \frac{\langle v + ry, b\rangle_L^2}{L[b]\, y}\Big), \quad \text{if } L[b] < 0. \tag{2.9}$$

If the integral representation for $H$ does not converge, we mean its analytic continuation with respect to $k$. If $D < 0$, one of the special functions above simplifies, and we find that

$$H(y; D) = \exp\big(-\pi D y / 2|L|\big)\, \Gamma\big(1 + \tfrac{l_\pm}{2} - k, -Dy/|L|\big),$$

where $\Gamma(s, y) = \int_y^\infty t^{s-1} e^{-t}\, dt$ is the upper incomplete gamma function.

Given $n \in \mathbb{R}$ and $r \in L_\mathbb{R}^\vee$, we call

$$D = D_L(n, r) := |L|\big(4n - L[r]\big)$$

the discriminant of $(n, r)$. We suppress dependence of $D_L(n, r)$ on $L$, $n$, and $r$, which will be clear from the context.

As an important side product to our study of Fourier expansion, we obtain further restrictions on non-trivial H-harmonic Maaß Jacobi forms.

**Proposition 2.11.** *Given an orthogonal set $B \subseteq \mathbb{P}(L_\mathbb{R})$ of non-isotropic vectors, then we have*

$$\mathcal{M}\mathbb{J}^{\mathrm{H}[B]}_{k,L} = \mathcal{M}\mathbb{J}^{\mathrm{h}[B_+]\mathrm{H}[B_-]}_{k,L},$$

*where $B_\pm = \{b \in B : \mathrm{sgn}(L[b]) = \pm 1\}$.*

**Proposition 2.12.** *Let $\phi \in \mathcal{M}\mathbb{J}^{\Delta \mathrm{h}[B_{\mathrm{hol}}] \mathrm{h}_0 \mathrm{H}[B]}_{k,L}$ be an H-harmonic Maaß-Jacobi with local Fourier expansion as in (2.6). With the same assumptions on $B_{\mathrm{hol}}$ and $B$ as above, for fixed $n$ and $r$, $c(\phi; n, r; y, v)$ it is an element of the complex vector space spanned by*

$$\prod_{b \in B'} H^{\mathrm{H}[b]}(y, v; r), \quad \text{and} \quad H(y; D) \prod_{b \in B'_-} H^{\mathrm{H}[b]}(y, v; r), \quad \text{for some } B' \subseteq B. \tag{2.10}$$

*Analogously, local Fourier expansions of $\phi \in \mathcal{M}\mathbb{J}^{\Delta \mathrm{h}[B_{\mathrm{hol}}] \overline{\mathrm{h}}_0 \mathrm{H}[B]}_{k,L}$ lie in the space spanned by*

$$e(2r(iv_0)) \prod_{b \in B'} H^{\mathrm{H}[b]}(y, v; r), \quad \text{and} \quad H(y; D) e(2r(iv_0)) \prod_{b \in B'_-} H^{\mathrm{H}[b]}(y, v; r), \quad \text{for some } B' \subseteq B. \tag{2.11}$$

The proof of Proposition 2.12 is a complicated calculation, parts of which we split up into three separate propositions.

**Proposition 2.13.** *Given $b \in L_\mathbb{R}$ that is is not isotropic, and a smooth function $a$, we have*

$$\Delta^{\mathrm{JH}[b]} a\left(\frac{\langle v, b\rangle_L^2}{4|L[b]|y}\right) = \frac{|L[b]| \Delta_{\frac{1}{2}} a(y) e(\mathrm{sgn}(L[b])\tau)}{y\, e(\mathrm{sgn}(L[b])\tau)}\bigg|_{y = \frac{\langle v,b\rangle_L^2}{4|L[b]|y}},$$

*where $\Delta_{1/2}$ is the classical weight $\frac{1}{2}$ Laplace operator defined in (1.10). Every solution to $\Delta^{\mathrm{JH}[b]}$ can be written as a function depending on $L$, $b$, $y$, and $v$ only by means of $\langle v, b\rangle_L^2 / 4|L[b]|y$.*

*More specifically, in the case of $L[b] < 0$ solutions to $\Delta^{\mathrm{JH}[b]} a(y, v; r)\zeta^r$ for some $r \in L_\mathbb{R}^\vee$ are linear combinations of the constant function and*

$$a(y, v; r) = \mathrm{sgn}\bigl(\langle v + r y, b\rangle_L\bigr) \gamma\left(\tfrac{1}{2}, -\pi \frac{\langle v + r y, b\rangle_L^2}{L[b] y}\right).$$

**Proposition 2.14.** *Let $B$ be as above, and fix $n \in \mathbb{R}$ and $r \in L_\mathbb{R}^\vee$. Given solutions $a_b(y, \langle v, b\rangle_L)$ to the differential equations*

$$\Delta^{\mathrm{JH}[b]} a_b(y, \langle v, b\rangle_L) \zeta^r = 0, \quad b \in B,$$

*set $a_B(y, v) = \prod_{b \in B} a_b(y, \langle v, b\rangle_L)$. Then the differential equations*

$$\mathscr{C}^{\mathrm{J}}_{k,L} a(y) a_B(y, v) q^n \zeta^r = 0 \quad \text{and} \quad \mathscr{C}^{\mathrm{J}}_{k,L} a(y) e(2r(iv_0)) a_B(y, v) q^n \zeta^r = 0$$

*are both equivalent to*

$$\Delta_{k - l_\pm/2} a(y) q^{D/4|L|} = 0.$$

**Proposition 2.15.** *Given $b \in L_\mathbb{R}$ that is not isotropic, let $a(y, v)$ be a real smooth function that is annihilated by $\Delta^{\mathrm{JH}[b]}$. Then $\xi^{\mathrm{JH}[b]}_L a(y, v)$ is constant, and $\xi^{\mathrm{J}}_{k,L} a(y, v) = 0$. More specifically, if $L[b] < 0$ then*

$$\xi^{\mathrm{JH}[b]}_L H^{\mathrm{H}[b]}(v, y; 0) = 2\sqrt{\pi}.$$

Proofs of Proposition 2.13, 2.14, and 2.15 can be found below.

*Proof of Proposition 2.11.* Assume that there is some $\phi \in \mathcal{M}\mathbb{J}_{k,L}^{\mathrm{H}[B]}$ and some $b \in B_+$ such that $\xi_L^{\mathrm{JH}[b]} \phi \neq 0$. We choose a maximal subset $\tilde{B}$ of $B$ that contains $b$ and satisfies

$$\xi_L^{\mathrm{JH}[\tilde{B}]} \phi \neq 0.$$

By Propositions 2.4 and 2.15, this is a non-trivial skew Jacobi form. Further, by Proposition 1.6, it cannot have any singularity in direction of $b$, which by a standard argument contradicts $L[b] > 0$. For details, the reader is refered to Theorem 1.2 in [EZ85]. It should be kept in mind that we use opposite signs for Jacobi indices. ∎

*Proof of Proposition 2.12.* Proposition 2.11 allows us to focus on the case $B = B_-$. Applying Proposition 2.13 to the case of negative $L[b]$ shows that any Fourier expansion is of the form

$$a(y) \prod_{b \in B_-} H^{\mathrm{H}[b]}(y, v; r).$$

Then Proposition 2.14 establishes the statement. ∎

*Proof of Proposition 2.13.* Recall that $\Delta^{\mathrm{JH}[b]} = -iy\big(i\partial_z(b) - 2\pi y^{-1}\langle\cdot, v\rangle_L\big)\partial_{\overline{z}}(b)$. A straightforward computation yields

$$\Delta^{\mathrm{JH}[b]} a\Big(\frac{\langle v,b\rangle_L^2}{4L[b]y}\Big) = \frac{\langle v,b\rangle_L^2}{4y} a''\Big(\frac{\langle v,b\rangle_L^2}{4L[b]y}\Big) + \Big(\frac{|L[b]|}{2} - \frac{\pi\langle v,b\rangle_L^2}{\mathrm{sgn}(L[b])y}\Big) a'\Big(\frac{\langle v,b\rangle_L^2}{4L[b]y}\Big).$$

On the other hand, $\Delta_{\frac{1}{2}} = 4y^2 \partial_\tau \partial_{\overline{\tau}} - iy\partial_{\overline{\tau}}$, from which we get

$$\Delta_{\frac{1}{2}} a(y) e\big(\mathrm{sgn}(L[b])\tau\big) = \Big(y^2 a''(y) + \big(\tfrac{y}{2} - 4\pi\mathrm{sgn}(L[b])y^2\big) a'(y)\Big) e\big(\mathrm{sgn}(L[b])\tau\big).$$

This completes the proof of the first part. The second part follows by combining covariance of $\Delta^{\mathrm{HJ}[b]}$ with respect to the action of $G^{(\mathrm{J}L)}$ and

$$\Delta_{\frac{1}{2}} \gamma\big(\tfrac{1}{2}, 4\pi y\big) e(-\tau) = 0.$$

∎

*Proof of Proposition 2.14.* Invariance of $\mathcal{C}_{k,L}^{\mathrm{J}}$ and $\Delta^{\mathrm{JH}[b]}$ under the action of the real Jacobi group, implies that we can restrict ourselfes to the case $r = 0$. Further, it suffices to treat the case of non-degenerate lattices $L$, since derivatives with respect to $L_0 \otimes \mathbb{R}$ do not occur in $\mathcal{C}^{\mathrm{J}}$. Then $D/4|L| = n$.

We first deduce three relations of differential operators, which hold under the mere assumption that $L[b] \neq 0$ for all $b \in B$. Sums and products in these formulas, if not indicated differently, run over elements in $B$. Since $L[b] \neq 0$, Proposition 2.13 implies that $a_b(y, v) = \tilde{a}_b\big(\langle v,b\rangle_L^2/2|L[b]|y\big)$ for some suitable function $\tilde{a}_b$. We deduce that $\partial_y a_b(y, v)$ equals

$$\frac{-\langle v,b\rangle_L^2}{4|L[b]|y^2} \tilde{a}_b'\Big(\frac{\langle v,b\rangle_L^2}{4|L[b]|y}\Big) = \frac{-\langle v,b\rangle_L}{4L[b]y} \frac{\langle v,b\rangle_L}{\mathrm{sgn}(L[b])y} \tilde{a}_b'\Big(\frac{\langle v,b\rangle_L^2}{4|L[b]|y}\Big) = \frac{-\langle v,b\rangle_L}{4L[b]y} \partial_v(b) a_b(y, v). \quad (2.12)$$

Applying Relation (2.12) to $\prod a_b(y, v)$ and using orthogonality of $B$, we find that

$$\partial_y \prod_{b \in B} a_b(y, v) = \sum_{b \in B} \frac{-\langle v, b \rangle_L \partial_v(b)}{4L[b]y} \prod_{\tilde{b} \in B} a_{\tilde{b}}(y, v) \tag{2.13}$$

$$\partial_y^2 \prod_{b \in B} a_b(y, v) = \Bigg( \sum_{b, b' \in B} \frac{\langle v, b \rangle_L \langle v, b' \rangle_L \partial_v(b) \partial_v(b')}{16 L[b] L[b'] y^2} + \sum_{b \in B} \frac{\langle v, b \rangle_L \partial_v(b)}{4 L[b] y^2} \Bigg) \prod_{\tilde{b} \in B} a_{\tilde{b}}(y, v). \tag{2.14}$$

In addition, using the equation $\Delta^{\mathrm{JH}[b]} a_b(y, v) = 0$, we infer that

$$\partial_v(b)^2 a_v(y, v) = \frac{4\pi \langle v, b \rangle_L \partial_v(b)}{y} a_b(y, v). \tag{2.15}$$

The actual proof of Proposition 2.14 is a computation, which we start by expanding Equation (2.14) for $\mathscr{C}^{\mathrm{J}}$ with respect to deriviatives $\partial_\tau$ and $\partial_{\overline{\tau}}$ that are applied to the factor $a(y) q^n$.

$$\mathscr{C}^{\mathrm{J}}_{k, L} a(y) q^n \prod_b a_b(y, v) = -2\big(\Delta_{k - l_\pm/2} a(y) q^n\big) \prod_b a_b(y, v)$$

$$- 8y^2 \big(\partial_\tau a(y) q^n\big) \partial_{\overline{\tau}} \prod_b a_b(y, v) \tag{2.16}$$

$$- 8y^2 \big(\partial_{\overline{\tau}} a(y) q^n\big) \partial_\tau \prod_b a_b(y, v) \tag{2.17}$$

$$+ \frac{y^2}{\pi i} \Big( \big(\partial_\tau a(y) q^n\big) L[\tfrac{i}{2} \partial_v] + \big(\partial_{\overline{\tau}} a(y) q^n\big) L[\tfrac{-i}{2} \partial_v] \Big) \prod_b a_b(y, v) \tag{2.18}$$

$$- 8y \big(\partial_\tau a(y) q^n\big) \langle v, \partial_{\overline{z}} \rangle_L \prod_b a_b(y, v) \tag{2.19}$$

$$+ a(y) q^n \mathscr{C}^{\mathrm{J}}_{k, L} \prod_b a_b(y, v). \tag{2.20}$$

In order to establish the proposition, it suffices to show that the sum of the last five terms (2.16) through (2.20) vanishes.

The sum of (2.16) and (2.17) equals

$$-2y^2 \big(2 a'(y) q^n - 4\pi n \, a(y) q^n\big) \partial_y \prod_{b \in B} a_b(y, v) = y \big(a'(y) - 2\pi n \, a(y)\big) q^n \sum_{b \in B} \frac{\langle v, b \rangle_L \partial_v(b)}{L[b]} \prod_{\tilde{b} \in B} a_{\tilde{b}}(y, v)$$

by Relation (2.13). In the term (2.18) terms involving $a'(y)$ cancel. To treat the derivative $L[\partial_v]$ applied to $\prod a_b(y, v)$, we recall one elementary fact. Since $B$ consists of mutually orthogonal vectors, we have

$$\partial_v \prod_{b \in B} a_b(y, v) = \sum_{b \in B} \frac{\langle b, \cdot \rangle_L}{2 L[b]} \partial_v(b) \prod_{\tilde{b} \in B} a_{\tilde{b}}(y, v). \tag{2.21}$$

Since $L[\langle b, \cdot \rangle_L] = 4 L[b]$, we find that (2.18) yields

$$\frac{-n y^2}{2} a(y) q^n L[\partial_v] \prod_{b \in B \in B} a_b(y, v) = \frac{-n y^2}{2} a(y) q^n \sum_{b \in B} \frac{4 L[b]}{4 L[b]^2} \partial_v(b)^2 \prod_{\tilde{b} \in B} a_{\tilde{b}}(y, v)$$

$$= -2\pi n y \, a(y) q^n \sum_{b \in B} \frac{\langle v, b \rangle_L \partial_v(b)}{L[b]} \prod_{\tilde{b} \in B} a_{\tilde{b}}(y, v)$$

by Relation (2.15). Finally, the term (2.19) can be expanded in a straightforward way:

$$-y\big(a'(y) - 4\pi n\, a(y)\big) q^n \sum_{b \in B} \frac{\langle v, b \rangle_L \partial_v(b)}{L[b]} \prod_{\tilde{b} \in B} a_{\tilde{b}}(y, v).$$

By what we have computed so far, we see that (2.16) + (2.17) + (2.18) + (2.19) = 0.

To prove the proposition we are therefore reduced to showing that (2.20) vanishes. This, in turn, is the same as showing that $\mathscr{C}^J_{k,L} \prod_b a_b(y, v) = 0$. Considering expression (1.8) for $\mathscr{C}^J_{k,L}$, we see that no terms involving $\partial_u$ contributes. The second term in (1.8) does not contribute either, because $\partial_\tau$ and $\partial_{\overline\tau}$ cancel each other. The order four term in (1.8) also does not contribute, as is easily seen by means of (2.21).

We are left with the following expression for $\mathscr{C}^J_{k,L} \prod a_b(y,v)$, which originates in the first, the third, the seventh, and the eighth term of (1.8).

$$\left( -2y^2 \partial_y^2 - (2k - l_\pm) y \partial_y - 2y \partial_y \langle v, \partial_v \rangle_L - \tfrac{1}{2} \langle v, \langle v, \partial_v \rangle_L, \partial_v \rangle_L - \tfrac{2k - l_\pm - 1}{2} \langle v, \partial_v \rangle_v \right) \prod_{b \in B} a_b(y,v).$$

We simplify it by employing (2.13), (2.14), and (2.21). It suffices to show vanishing of:

$$-\sum_{b,b'} \frac{\langle v,b \rangle_L \langle v,b' \rangle_L \partial_v(v)\, \partial_v(b')}{8 L[b] L[b']} - \sum_{b \in B} \frac{\langle v,b \rangle_L \partial_v(b)}{2 L[b]} + (2k - l_\pm) \sum_{b \in B} \frac{\langle v,b \rangle_L \partial_v(b)}{4 L[b]}$$

$$+ \sum_{b,b' \in B} \frac{\langle b', v \rangle_L \partial_v(b')}{2 L[b']} \frac{\langle b, v \rangle_L \partial_v(b)}{2 L[b]} - \frac{1}{2} \sum_{b,b' \in B} \frac{\langle b,v \rangle_L \langle b',v \rangle_L \partial_v(b) \partial_v(b')}{4 L[b] L[b']} - \frac{2k - l_\pm - 1}{2} \sum_{b \in B} \frac{\langle v,b \rangle_L \partial_v(b)}{2 L[b]}$$

The first, the order 2 contribution of the fourth, and the fifth term cancel. The order 1 contribution of the fourth term is $\sum_b \langle b, v \rangle_L \partial_v(b)/2L[b]$, so that the remaining terms cancel, too. This completes the proof. ∎

*Proof of Proposition 2.15.* The statement on $\xi^J$ is a consequence of (2.21) in the proof of Proposition 2.14.

Our proof of the first statement on $\xi^{JH[b]}$ uses Proposition 2.13. Let $a(y)$ be a solution to the differential equation $\Delta_{\frac{1}{2}} a(y) e\big(\mathrm{sgn}(L[b])\tau\big)$. If $a(y) \ne 1$, we can choose a scalar multiple of $a(y)$ in such a way that

$$\xi_{\frac{1}{2}} a(y) e\big(\mathrm{sgn}(L[b])\tau\big) = e\big(-\mathrm{sgn}(L[b])\tau\big), \quad \text{therefore} \quad a'(y) = y^{\frac{-1}{2}} \exp\big(4\pi\,\mathrm{sgn}(b)\, y\big).$$

The proposition then follows from a straightforward computation. ∎

# 3 Indefinite theta series

In this section, we slightly generalize the results of Chapter 2 of [Zwe02]. Zwegers defined indefinite theta series for lattices of signature $(r+1, 1)$, $r \ge 0$. These theta series are H-harmonic, as can be deduced from Proposition 2.12. Zwegers also remarked that his construction generalizes to arbitrary, non-degenerate lattices. Here, we make a small step towards the most general type of indefinite theta series by considering the ones of product type (this terminology is explained in the introduction). In the previous section, we have seen that in our setting we may only expect real analytic contributions from the negative definite part of $L$. Further, such real analytic terms in the H-harmonic setting must be associated to mutually orthogonal directions in $L$. One possibility to achieve this is to consider contribution from orthogonal subspaces of $L_\mathbb{R}$ that have signature $(1,1)$. Many cases that we consider here could be alternatively constructed using products of Zwegers's indefinite theta series. However, in Section 3.1 we show that this is not always true. On the other hand, we would like to emphasize that the intention behind this section is not to construct all indefinite theta series. In fact, even when restricting to H-harmonic indefinite theta series, the presented construction does not exhaust all cases. Rather, in the present section we provide some non-trivial examples to bring the previously developed theory to life.

Throughout this section, we assume that $L$ is non-degenerate. We repeat Zwegers's construction on orthogonal pieces of signature $(1,1)$, given by a set

$$C = \{(c_1, c_2) \in \mathbb{P}(L_{\mathbb{R}})^2 : \mathrm{sgn}(\mathrm{span}\{c_1, c_2\}) = (1,1), \langle c_1, c_2 \rangle_L < 0\}$$

of size $l_-$. For any distinct pairs $(c_1, c_2)$ and $(c'_1, c'_2)$ in $C$, we assume that $\mathrm{span}\{c_1, c_2\}$ is orthogonal to $\mathrm{span}\{c'_1, c'_2\}$. A set $C$ with this property we be called a $(1,1)$-decomposition of $L$. As an obvious restriction on the existence of $(1,1)$-decompositions, we have $l_+ \geq l_-$. Indeed, $\mathrm{span}\, C \subseteq L_{\mathbb{R}}$ has signature $(l_-, l_-)$.

To define indefinite theta series, we need an analogue of Zwegers's $\rho$ function for $c \in \mathbb{P}(L_{\mathbb{R}})$ with $L[c] \leq 0$:

$$\begin{aligned}\rho^c(\tau, z; \nu) &= \mathrm{sgn}(\langle c, \nu \rangle_L), && \text{if } L[c] = 0; \text{ and} \\ \rho^c(\tau, z; \nu) &= H^{\mathrm{H}[c]}(y, v; \nu), && \text{if } L[c] < 0.\end{aligned} \quad (3.1)$$

**Definition 3.1.** Let $C$ be a $(1,1)$-decomposition of $L$. The vector valued (indefinite) theta series of $L$ attached to $C$ is defined by

$$\theta_L^C(\tau, z) = \sum_{\overline{\nu} \in \mathrm{disc}\, L} \mathfrak{e}_{\overline{\nu}} \sum_{\nu \in \overline{\nu} + L} e(L[\nu]\tau + \langle z, \nu \rangle_L) \prod_{(c_1, c_2) \in C} (\rho^{c_1} - \rho^{c_2})(\tau, z; \nu). \quad (3.2)$$

*Example 3.2 (Zwegers's indefinite theta series for lattices of signature $(r-1,1)$).* Zwegers analyzed indefinite theta series for Lorentzian lattices. For the time being, we adopt notation from Chapter 2 of [Zwe02]. Starting with a quadratic space $V$ of signature $(r-1,1)$, $r \geq 2$, choose isotropic or negative vectors $c_1, c_2 \in \mathbb{R}^r$. For any such choice, one obtains an indefinite theta series

$$\theta^{c_1, c_2}(\tau, z) = \sum_{\nu \in \mathbb{Z}^r} \rho^{c_1, c_2}(\nu; \tau) e(L[\nu]\tau + \langle z, \nu \rangle_L),$$

where, in the present paper's notation, we have

$$\rho^{c_1, c_2}(\nu; \tau) = \rho^{c_1}(\tau, z; \nu) - \rho^{c_2}(\tau, z; \nu).$$

We can recover $\theta^{c_1, c_2}$ as the $\mathfrak{e}_0$-component of $\theta_L^C$ with $C = \{(c_1, c_2)\}$ as long as $c_1$ is not a multiple of $c_2$. If $c_1 = c_2$ as elements of $\mathbb{P}(V_{\mathbb{R}})$, then $\theta^{c_1, c_2} = 0$.

Indefinite theta series have singularities, whose location is determined by the geometry of $C$. We set

$$D(C) := \{(\tau, z) \in \mathbb{H}^{(JL)} : \langle b_1, v/y \rangle_L, \langle b_2, v/y \rangle_L \notin \mathbb{Z} \text{ for all } (b_1, b_2) \in C\}.$$

**Proposition 3.3.** *Given a $(1,1)$-decomposition of $L$, the theta series $\theta_L^C$ converges absolutely and locally uniformly for any $z \in D(C)$. It can be analytically continued to $\mathbb{H}^{(JL)}$ except for non-moving singularities of almost meromorphic type.*

The proof of Proposition 3.3 relies on a reduction to rational $c_1$'s and $c_2$'s. The next lemma allows us, in addition, to replace isotropic vectors by negative ones.

**Lemma 3.4.** *Let $C$ be a $(1,1)$-decomposition of $L$. Fix $(c_1, c_2) \in C$ such that for $i = 1$ and $i = 2$ we have $\mathbb{R}c_i \cap L_\mathbb{Q} = \emptyset$ or $L[c_i] = 0$. Then there is $\tilde{c} \in \mathbb{P}(L_\mathbb{Q})$ with $L[\tilde{c}] < 0$ such that both*

$$(\tilde{c}, c_2) \cup C \setminus \{(c_1, c_2)\} \quad \text{and} \quad (c_1, \tilde{c}) \cup C \setminus \{(c_1, c_2)\}$$

*are $(1,1)$-decompositions of $L$.*

*Proof.* We analyze the following Grassmannian of negative definite, rational subspaces

$$H = \mathrm{Gr}_- \{v \in L_\mathbb{Q} : \forall (c_1', c_2') \in C, (c_1', c_2') \neq (c_1, c_2) : \langle v, c_i' \rangle_L = 0\}.$$

Since $\mathrm{span}\{c_1, c_2\}$ has signature $(1,1)$, $(C \setminus \{c_1, c_2\})^\perp$ has signature $(l_+ - l_- + 1, 1)$. Therefore $H$ is a manifold of dimension at least 1. By our assumption on $c_1$ and $c_2$, we find that $\mathbb{R}c_1, \mathbb{R}c_2 \notin H$. It thus suffices to choose any $\mathbb{Q}\tilde{c} \in H$. ∎

The next proof is essentially due to Zwegers [Zwe02]. A little care must be taken when splitting off $\mathrm{span}\{c_1, c_2\}$ for $(c_1, c_2) \in C$, and for this reason we give some details. For our purpose, it is also important that we prove that $\theta_L^C$ has non-moving singularities of almost meromorphic type. This fact was not mentioned in [Zwe02], even though it is immediate from the treatment given there.

*Proof of Proposition 3.3.* Using Lemma 3.4, we can successively write $\theta_L^C$ as the sum of $\theta_L^{\widetilde{C}}$'s where for each $(c_1, c_2) \in \widetilde{C}$ at least one of the $c_i$'s is rational and negative. After swapping $c_1$ and $c_2$ if necessary, we can thus assume without loss of generality that $L[c_1] < 0$ and $c_1 \in L_\mathbb{Q}$ for all $(c_1, c_2) \in C$.

It also suffices to show convergence for $\overline{\nu} = 0$. Indeed, given $\overline{\nu} \in \mathrm{disc}(L)$, we can replace $z$ by $z + \nu$ for a representative $\nu$ of $\overline{\nu}$. The set $D(C)$ remains unchanged under this substitution, since $\nu \in L^\vee$.

Under these hypotheses, we can proceed in a similar way as Zwegers. It will be convenient to write $z$ as $a + \tau b$ for $a, b \in L_\mathbb{R}$. In particular, we have $\nu/y = a$. Let $\beta(y) = \int_y^\infty t^{\frac{-1}{2}} e^{-\pi t} \, dt$ be the beta function. We have

$$H^{\mathrm{H}[c]}(y, \nu; \nu) = \sqrt{\pi} \, \mathrm{sgn}(\langle a + \nu, c \rangle_L) \left(1 - \beta\big(-y \langle a + \nu, c \rangle_L^2 / L[c]\big)\right).$$

So, up to scalar multiples of $\sqrt{\pi}$, we can write $\rho^{c_1} - \rho^{c_2}$ as a sum or difference of terms

$$\mathrm{sgn}(\langle a+\nu, c_1 \rangle_L) \beta\big(-y\langle a+\nu, c_1\rangle_L^2 / L[c_1]\big), \quad \mathrm{sgn}(\langle a+\nu, c_2\rangle_L) \beta\big(-y\langle a+\nu, c_2\rangle_L^2 / L[c_2]\big), \quad \text{and} \tag{3.3}$$

$$\mathrm{sgn}(\langle a+\nu, c_1\rangle_L) - \mathrm{sgn}(\langle a+\nu, c_2\rangle_L).$$

We have to estimate the product $\prod_{(c_1, c_2)} (\rho^{c_1} - \rho^{c_2})$. The orthogonality relation imposed on elements of $C$ allows us to focus on each term individually, so that Zwegers's estimates apply word by word for negative $c_1$ or $c_2$.

We can thus focus on the case that $c_1 \in L_\mathbb{Q}$ is negative and $c_2$ is isotropic. For simplicity write $C_0 \subseteq C$ for the set $(c_1, c_2)$ with isotropic $c_2$. We discuss the contribution of the third term in (3.3) in more detail, since we need it to understand singularities of $\theta_L^C$. We are free to replace $c_1$ by a scalar multiple of itself without changing $\theta_L^C$. In particular, we may assume, by Lemma 3.4, that $c_1 \in L$.

Given $(c_1, c_2) \in C_0$, we decompose $a + \nu = \mu + nc_1$ with $n \in \mathbb{Z}$ and $0 \le \langle c_2, \mu \rangle_L < \langle c_2, c_1 \rangle_L$. As in the proof of Proposition 2.4 in [Zwe02], we can use the following equality, since $z \in D(C)$.

$$\sum_{\substack{n \in \mathbb{Z} \\ (c_1, c_2) \in C}} \left( \mathrm{sgn}(\langle \mu, c_2 \rangle_L) - \mathrm{sgn}\left(n + \frac{\langle \mu, c_1 \rangle_L}{\langle c_1, c_2 \rangle_L}\right) \right) e(\langle \mu, c_2 \rangle n\tau + \langle b, c_2 \rangle n)$$
$$= \frac{2}{1 - e(\langle \mu, c_2 \rangle_L \tau + \langle b, c_2 \rangle_L)} - \delta(\langle \mu, c_1 \rangle_L), \quad (3.4)$$

where $\delta(0) = 1$ and $\delta(t) = 0$ for $t \neq 0$.

Because the $(c_1, c_2) \in C$ are mutually orthogonal, we can in fact decompose $a + \nu = \mu + \sum_{(c_1, c_2) \in C_0} n_{c_1}$. By passing to the corresponding quotient of $L$, we see that

$$\mathcal{M}(C_0) = \left\{ \mu \in a + L : \forall (c_1, c_2) \in C_0 : 0 \le \langle c_2, \mu \rangle_L < \langle c_2, c_1 \rangle_L \right\} / \left( \mathrm{span}\{c_1 : (c_1, c_2) \in C_0\} \right)^\perp$$

is finite. Up to sign, the contribution to $\theta_L^C$ under consideration thus equals

$$\sum_{\substack{\mu \in \mathcal{M}(C_0) \\ \nu \in (\mathrm{span}\{c_1 : (c_1, c_2) \in C_0\})^\perp}} e(L[\mu + \nu]\tau + \langle z, \mu + \nu \rangle_L) \prod_{(c_1, c_2) \in C \setminus C_0} (\rho^{c_1} - \rho^{c_2})(\tau, z, \mu + \nu)$$
$$\cdot \prod_{(c_1, c_2) \in C_0} \left( \frac{2}{1 - e(\langle \mu + \nu, c_2 \rangle_L \tau + \langle b, c_2 \rangle_L)} - \delta(\langle \mu, c_1 \rangle_L) \right). \quad (3.5)$$

This reduces us to convergence of $\theta_{L'}^{C \setminus C_0}$ with $L' = \left( \mathrm{span}\{c_1 : (c_1, c_2) \in C_0\} \right)^\perp$, for which we have already referred the reader to [Zwe02]. From this point on, it is clear how to establish convergence of $\theta_L^C$ as in [Zwe02].

We are left with proving that $\theta_L^C$ has non-moving singularities of almost meromorphic type. Observe that (3.4) can be continued to all except finally many $b \pmod{L}$. Therefore, singularities of (3.5) arise from products of meromorphic functions with real-analytic ones. This completes the proof. ∎

We can go one step beyond Proposition 3.3, by establishing uniform convergence. The next lemma will help us proving Theorem 3.6.

**Lemma 3.5.** *Fix a $(1, 1)$-decomposition $C$ of $L$, and one pair $(c_1, c_2) \in C$. Let $\tilde{c}$ be a continuous function $\tilde{c} : [0, \infty) \to \mathrm{span}\{c_1, c_2\}$. Assume that $L[\tilde{c}(t)] < 0$, $L[c_2] = 0$, and $\tilde{c}(0) = c_2$, and that $\mathrm{span}\{c_1 \tilde{c}(t)\}$ has signature $(1, 1)$ for all $t > 0$.*

*Set $C(t) = \{(\tilde{c}(t), c_2)\} \cup C \setminus \{(c_1, c_2)\}$. Then $\lim_{t \to 0} \theta_L^{C(t)} = \theta_L^C$.*

*Proof.* The proof of Proposition 2.7 (4) in [Zwe02] applies almost word by word. We give some details for convenience.

In analogy to the proof of Proposition 3.3, we can assume that $\overline{\nu} = 0$. Zwegers's argument saying that we can restrict to the situation $c_1 = c_2$ and $\mathrm{span}\{c_1, \tilde{c}(t)\}$ be constant for all $t > 0$ $\tilde{c}(t)$ can be transferred to our setting without difficulty. We thus have to show that

$$\sum_{\nu \in \overline{\nu} + L} e(L[\nu]\tau + \langle z, \nu \rangle_L) \prod_{(c_1', c_2') \in C(t)} (\rho^{c_1'} - \rho^{c_2'})(\tau, z; \nu) \longrightarrow 0 \quad \text{as } t \to 0.$$

Writing $\rho^{c_1} - \rho^{\tilde{c}(t)}$ as a sum of

$$\mathrm{sgn}(\langle a + \nu, \tilde{c}(t) \rangle_L) \beta\left( -y \langle a + \nu, \tilde{c}(t) \rangle_L^2 / L[\tilde{c}(t)] \right) \quad \text{and}$$
$$\mathrm{sgn}(\langle a + \nu, c_1 \rangle_L) - \mathrm{sgn}(\langle a + \nu, \tilde{c}(t) \rangle_L),$$

we are reduced to two separate cases. It suffices to apply the inequalities

$$\left|\operatorname{sgn}(\langle a+v,c_1\rangle_L) - \operatorname{sgn}(\langle a+v,c_2\rangle_L)\right| \le \left|\operatorname{sgn}(\langle a+v,c_1\rangle_L) - \operatorname{sgn}(\langle a+v,\tilde{c}(t)\rangle_L)\right|,$$

and

$$\left|\beta\bigl(-y\langle a+v,\tilde{c}(t)\rangle_L^2/L[\tilde{c}(t)]\bigr)e\bigl(L[a+v]\tau + \langle a+v,z\rangle_L\bigr)\right| \le e\bigl(\widetilde{L}[a+v]\tau + \langle a+v,z\rangle_{\widetilde{L}}\bigr)$$

for a suitable quadratic form $\widetilde{L}$, which depends on $\langle c_1, a+v\rangle_L$ and $\langle c_2, a+v\rangle_L$. Both estimates were literally established by Zwegers. We close by revisiting the construction of $\widetilde{L}$.

Zwegers splits up the set $L$ into three smaller ones, which he calls $P_1$, $P_2$, $P_3$. One $P_1$ and $P_2$ one quickly obtains a majorant on $P_3$, however, he needs to use additional vectors $\tilde{c}_1 = \frac{\langle c_1, c_2\rangle_L}{2L[c_1]} c_1 - c_2$ and $\tilde{c}_2 = -c_2$ and a quadratic form

$$\widetilde{L}[v] = L[v^\perp] - \frac{\langle c_2, v\rangle_L\bigl(\langle c_1, c_2\rangle_L \langle c_1, v\rangle_L - L[c_1]\langle c_2, v\rangle_L\bigr)}{\langle c_1, c_2\rangle_L^2}.$$

The concrete estimate in the last case is involved, but we nevertheless urge the reader to consult [Zwe02] on this matter. ∎

If $L$ is an integral lattice, then we immediately see that $\theta_L^C$ is invariant under the action of the Heisenberg group in $\Gamma^{(JL)}$. The correct transformation behavior with respect to $\operatorname{SL}_2(\mathbb{Z})$ can be inferred, again, in the same way as in [Zwe02]. Zwegers, however, reminded the author of Vignéras's results [Vig77], which we achieved to write in a particularly neat way.

**Theorem 3.6.** *For every $L$ and every $C$, the theta series $\theta_L^C$ transforms like a Jacobi form of weight $l/2$, index $L$, and type $\rho_L^J$ with respect to the action of $\operatorname{SL}_2(\mathbb{Z}) \subset \Gamma^{(JL)}$.*

**Corollary 3.7.** *Suppose that $L$ is an even lattice. Then the the theta series $\theta_L^C$ transforms like a Jacobi form of weight $l/2$, index $L$, and type $\rho_L^J$. If $L[c_2] = 0$ for all $(c_1,c_2) \in C$, then*

$$\theta_L^C \in \mathcal{M}\mathbb{J}_{\frac{l}{2},L}^{\mathrm{H}[B]}(\rho_L), \quad \text{where} \quad B = \{c_1 : (c_1,c_2) \in C\}.$$

*Proof of Theorem 3.6.* In light of Lemma 3.5, we can reduce ourselves to the case that $c_1$ and $c_2$ are negative for all $(c_1,c_2) \in C$. The estimates in Proposition 3.3, i.e. in Proposition 2.4 of [Zwe02], show that

$$e\bigl(L[v]\tau + \langle z,v\rangle_L\bigr) \prod_{(c_1,c_2)\in C} \bigl(\rho^{c_1} - \rho^{c_2}\bigr)(\tau,z;v)$$

is of exponential decay. In particular, it is integrable and square integrable with respect to $v$, and so are its differentials and products with arbitrary polynomials in $v$. Vignéras states that $\theta_L^C$ is modular, if the above term at $y = 1$ is annihilated by $\frac{1}{4\pi}\Delta_L - E_L$, where $\Delta_L$ is the Laplace operator attached to $L$ and $E_L$ is the corresponding Euler operator. One verifies that this, up to scalar constants, equals $\langle \mathrm{R}^{\mathrm{JH}}, \mathrm{L}^{\mathrm{JH}}\rangle_L$ restricted to $y = 1$.

It suffices to verify the differential equation after expanding the product over $C$. So fix $i_j \in \{1,2\}$ for $1 \le j \le l_-$ and set $B = \{c_{i_j} : (c_1,c_2) \in C\}$. Then

$$\langle \mathrm{R}^{\mathrm{JH}}, \mathrm{L}^{\mathrm{JH}}\rangle_L = \langle \mathrm{R}^{\mathrm{JH}}, \mathrm{L}^{\mathrm{JH}}\rangle_{(\operatorname{span} B)^\perp} + \sum_{b \in B} \Delta^{\mathrm{JH}[b]}.$$

By construction $\prod_{b \in B} \rho^b$ has H-harmonicity $B$, and this completes the proof. ∎

**§ 3.1 Reduction to smaller lattices.** The relation between indefinite theta series of product type and products of indefinite theta series is not quite obvious. We give a sufficient criterion that asserts splitting into more than one factor. It is crucial to note that in the assumptions of the next proposition we intersect the span of $C_j \subset C$ with a *rational* subspace of $L$.

**Proposition 3.8.** *Let $L$ and $C$ be as above. Suppose that $C$ can be written as a disjoint union of $C_j$, with $1 \le j \le l_C$. Suppose that further $L_\mathbb{Q} \supset \bigoplus_{j=1}^{l_C} L_j$ for lattices $L_j$ such that $l = \sum_j l_j$. If $\mathrm{span}\, C_j \cap L_{j,\mathbb{Q}} = L_{j,\mathbb{Q}}$, then*

$$\pi_L\big(\bigoplus \theta_{L_j}^{C_j}\big) = \theta_L^C,$$

*where $\pi_L$ is the natural projection of $\bigotimes_j \rho_{L_j}$ onto $\rho_L$ that is induced by the inclusion $\mathrm{disc}(L) \subseteq \bigotimes \mathrm{disc}(L_j)$.*

*Proof.* It is clear how to reduce ourselves to the case that $L = \bigoplus_j L_j$ by means of the projection $\pi_L$. Then it suffices inspect the definition (3.2) of $\theta_L^C$, to prove the proposition. ∎

We give two examples of indefinite theta series, which illustrate extreme cases that can occur. Both are attached the same lattice of signature $(2,2)$, but we use different $C$. The first one is a product of two indefinite theta series, the second one cannot possibly split.

*Example 3.9.* Let $L_0$ be the lattice with Gram matrix $\begin{pmatrix} 2 & 0 \\ 0 & -2 \end{pmatrix}$, and set $L = L_0 \oplus L_0$. To simplify notation, we further fix corresponding bases and coordinates $z_1, z_2$ and $z_1, z_2, z_3, z_4$ for $z$, respectively.

(1) Set $c_{0,1} = (-1,1)$, $c_{0,2} = (0,1)$, and $C_0 = \{(c_{0,1}, c_{0,2})\}$. We consider the product of two of Zwegers's indefinite theta series attached to $L_0$ and $C_0$. It yields an H-harmonic Maaß-Jacobi form:

$$\theta_L^C(\tau, z) = \theta_{L_0}^{C_0}(\tau, z_1, z_2) \cdot \theta_{L_0}^{C_0}(\tau, z_3, z_4) \quad \text{with} \quad C = \{\big((-1,1,0,0),(0,1,0,0)\big), \big((0,0,-1,1),(0,0,0,1)\big)\}.$$

Its H-harmonicity is $B = \{(0,1,0,0),(0,0,0,1)\}$, and it is a product of two indefinite theta series.

(2) The same lattice $L$ allows for an indefinite theta series that is not a product (but per terminology is of product type).

$$\theta_L^C(\tau, z) \quad \text{with} \quad C = \{\big((1,1,\pi^2,\pi^2),(1,\pi^2,1,\pi)\big), \big((\pi^2,\pi^2,1,1),(1,\pi,1,\pi^2)\big)\}.$$

This theta series is H-harmonic, because $(1,1,\pi^2,\pi^2)$ and $(\pi^2,\pi^2,1,1)$ are isotropic. To check that it cannot split as a product, it suffices to observe that

$$\mathrm{span}_{\mathbb{R}}\{(1,1,\pi^2,\pi^2),(1,\pi^2,1,\pi)\} \cap L_\mathbb{Q}$$

has dimension 0, while the vectors in $C$ span $L_\mathbb{R}$.

*Remark 3.10.* The preceding example makes clear that one has to carefully distinguish between product of indefinite theta series and indefinite theta series of product type. The latter terminology can be explained by alluding to a jet to be developed adelic theory of mock theta series. At $\infty$ such a theory should not distinguish between the two concepts, and the different that was showcased in Example 3.9 should originate in the finite places. This also explains why there is such a huge overlap between the proofs here and those in [Zwe02]: Convergence and even modularity of theta series, is mostly an issue that can be handled "at the infinite places".

## 4 Theta-like decompositions

We conclude this work with a theta-like decomposition for H-harmonic Maaß-Jacobi forms. Generalized $\widehat{\mu}$ functions play the role of $\theta$ functions.

**§4.1 Generalized $\widehat{\mu}$-functions.**  Indefinite theta series can be used to obtain preimages of skew theta series under $\xi_L^{\mathrm{JH}[B]}$. This construction, in the case of $l_- = 1$, should be somewhat known to experts. Write $\langle \cdot, \mathfrak{e}_{\overline{\nu}} \rangle_\mathfrak{e}$ for the $\overline{\nu}$-th coordinate of a vector valued function. Given a negative definite, even lattice $L$ and an orthogonal basis $B$ of $L$, we define

$$\widehat{\mu}_L^B(\tau,z) = \frac{\sum_{\overline{\nu}\in\mathrm{disc}\,L}\mathfrak{e}_{\overline{\nu}}\langle\theta^C_{L\oplus L(-1)}(\tau,(z,0)),\mathfrak{e}_{\overline{\nu}}\mathfrak{e}_0\rangle_\mathfrak{e}}{\langle\theta_{L(-1)}(\tau,0),\mathfrak{e}_0\rangle_\mathfrak{e}}, \tag{4.1}$$

where $C = \{((b,0),(b,b)) : b \in B\}$. Clearly $(0,b)$ is negative definite, while $(b,b)$ is isotropic. Singularities of $\theta^C_{L\oplus L(-1)}(\tau,(z,0))$ lie outside of $D(B) = \{(\tau,z) \in \mathbb{H}^{(\mathrm{J}L)} ; \langle b, \nu/y \rangle_L \notin \mathbb{Z}$ for all $b \in B\}$.

**Proposition 4.1.**  *Given $L$ and $B$ as above, we have*

$$\xi_L^{\mathrm{JH}[B]} \widehat{\mu}_L^B = 2\sqrt{\pi}\,\theta_L^{\mathrm{span}\,B}.$$

*Proof.* It suffices to determine the image of $\rho^c$ in (3.1) under $\xi_L^{\mathrm{JH}[b]}$. This was done in Proposition 2.15. ∎

**Theorem 4.2.**  *Suppose that $L = L' \oplus \bigoplus_{b\in B}\mathbb{Z}b$. Then we have*

$$\mathscr{M}\mathbb{J}_{k,L}^{\mathrm{H}[B]} = \sum_{B'\subseteq B} \left\langle \widehat{\mu}_L^{B'}, \mathscr{M}\mathbb{J}_{k-\#B'/2,(\mathrm{span}\,B')^\perp}\big(\rho_{\mathrm{span}\,B'}\big)\right\rangle. \tag{4.2}$$

*Remark 4.3.* (1) The condition that $L$ be $L' \oplus \bigoplus \mathbb{Z}b$ might seem surprisingly restrictive. However, for every lattice $L$ and ever $B \subseteq \mathbb{P}(L_\mathbb{Q})$ we can pass to such a situation by considering suitable suplattices $L_\perp \supseteq L$ in conjunction with the theory of vector valued Jacobi forms. On the other hand, overlapping singularities of meromorphic Jacobi forms in decomposition (4.2) hinder us from passing back from $L_\perp$ to $L$.

(2) A similar decomposition can be found for $\mathscr{M}\mathbb{J}_{k,L}^{\mathrm{h_0 H}[B]}$ and $\mathscr{M}\mathbb{J}_{k,L}^{\overline{\mathrm{h}}_0\mathrm{H}[B]}$, in the case of degenerate $L$.

*Proof of Theorem 4.2.* The theorem follows by induction: When applying $\xi^{\mathrm{JH}[B]}$ to $\phi \in \mathscr{M}\mathbb{J}_{k,L}^{\mathrm{H}[B]}$ to obtain a vector valued meromorphic Jacobi form by means of the partial theta decomposition Proposition 1.7. Using $\widehat{\mu}_L^B$ we can find an H-harmonic Maaß-Jacobi form $\psi$ that has the same image under $\xi^{\mathrm{JH}[B]}$ as $\phi$. Thus the difference $\phi - \psi$ is a sum of H-harmonic Maaß-Jacobi forms with H-harmonicities $B'$, $\#B' < \#B$. This completes the proof. ∎

**§4.2 Restriction to torsion points.**  To formulate the final corollary, we have to define harmonic weak Maaß forms of higher depth. It does not seem adequate to employ the generally preferable language of vector valued modular forms. Technicalities would require an extra exposition. Therefore, given $\Gamma \subseteq \mathrm{SL}_2(\mathbb{Z})$, we set $\mathbb{M}_k^{[0]}(\Gamma) = \mathrm{M}_k(\Gamma)$. For depth $d \geq 1$, we let $\mathbb{M}_k^{[d]}(\Gamma)$ be the space of real-analytic functions $f : \mathbb{H} \to \mathbb{C}$ that

(i) are invariant under the weight $k$ slash action of $\Gamma$;

(ii) whose image under $\xi_k$ lies in $\sum_l \overline{\mathrm{M}_l(\Gamma)} \otimes \mathbb{M}_{k-l}^{[d-1]}(\Gamma)$;

(iii) satisfy the growth condition $f(iy) = O(e^{ay})$ as $y \to \infty$ for some $a \in \mathbb{R}$.

The definition of higher depth harmonic weak Maaß forms goes back to Zagier and Zwegers, but currently there is no literature on it.

**Corollary 4.4.** *Suppose that* $(\tau, \alpha\tau + \beta)$ *for* $\alpha, \beta \in \frac{1}{N}L$ *is not a singularity of* $\phi \in \mathcal{MJ}_{k,L}^{H[B]}$. *Then we have*

$$\left(\phi|_{k,L}\left(\begin{pmatrix} 1 & 0 \\ 0 & 1 \end{pmatrix}, \alpha, \beta\right)\right)(\tau, 0) \in \mathbb{M}_k^{[\#B]}(\Gamma(N)),$$

*where* $\Gamma(N)$ *is the principal congruence subgroup of level* $N$.

*Remark 4.5.* The statement could be refined even more by specifying that each time we apply $\xi$ to the restriction of an H-harmonic Maaß-Jacobi form, a unary theta series splits off.

*Proof.* The transformation behavior is clear. The analytic properties follow from restricting to torsion points the Fourier expansions in Proposition 2.6. ∎

# References

[BCR14]   K. Bringmann, T. Creutzig, and L. Rolen. *Res. Math. Sci.* 1 (2014).

[BF04]    J. H. Bruinier and J. Funke. "On two geometric theta lifts". *Duke Math. J.* 125.1 (2004).

[BMR14]   K. Bringmann, K. Mahlburg, and R. C. Rhoades. "Taylor coefficients of mock-Jacobi forms and moments of partition statistics". *Math. Proc. Cambridge Philos. Soc.* 157.2 (2014).

[BO15]    K. Bringmann and R. Olivetto. *Kac-Wakimoto Characters and Non-Holomorphic Jacobi Forms.* Bringmann's homepage. To appear in Trans. Amer. Math. Soc. 2015.

[BR10]    K. Bringmann and O. K. Richter. "Zagier-type dualities and lifting maps for harmonic Maass-Jacobi forms". *Adv. Math.* 225.4 (2010).

[BRR12a]  K. Bringmann, M. Raum, and O. K. Richter. *Harmonic Maass-Jacobi forms with singularities and a theta-like decomposition.* arXiv:1207.5600. 2012.

[BRR12b]  K. Bringmann, M. Raum, and O. K. Richter. "Kohnen's limit process for real-analytic Siegel modular forms". *Adv. Math.* 231.2 (2012).

[Bru02]   J. H. Bruinier. *Borcherds products on O(2, l) and Chern classes of Heegner divisors.* Vol. 1780. Lecture Notes in Mathematics. Berlin: Springer-Verlag, 2002.

[BRZ15]   K. Bringmann, L. Rolen, and S. Zwegers. *On the Fourier coefficients of negative index meromorphic Jacobi forms.* arXiv:1501.04476. 2015.

[BS98]    R. Berndt and R. Schmidt. *Elements of the representation theory of the Jacobi group.* Vol. 163. Progress in Mathematics. Basel: Birkhäuser Verlag, 1998.

[CR10]    C. Conley and M. Raum. *Harmonic Maaß–Jacobi forms of degree* 1 *with higher rank indices.* arXiv:1012.2897. 2010.

[DIT11]   W. Duke, Ö. Imamoğlu, and Á. Tóth. "Cycle integrals of the $j$-function and mock modular forms". *Ann. of Math. (2)* 173.2 (2011).

[DMZ12]   A. Dabholkar, S. Murthy, and D. B. Zagier. *Quantum Black Holes, Wall Crossing, and Mock Modular Forms.* arXiv: 1208.4074. 2012.

[EZ85]    M. Eichler and D. B. Zagier. *The theory of Jacobi forms.* Vol. 55. Progress in Mathematics. Boston, MA: Birkhäuser Boston Inc., 1985.

[GZ98]    L. Göttsche and D. B. Zagier. "Jacobi forms and the structure of Donaldson invariants for 4-manifolds with $b_+ = 1$". *Selecta Math. (N.S.)* 4.1 (1998).

[Hel59]   S. Helgason. "Differential operators on homogenous spaces". *Acta Math.* 102 (1959).

[HL15]    K. Hikami and J. Lovejoy. "Torus knots and quantum modular forms". *Res. Math. Sci.* 2 (2015).

[Koh94]   W. Kohnen. "Nonholomorphic Poincaré-type series on Jacobi groups". *J. Number Theory* 46.1 (1994).

[KW01]   V. G. Kac and M. Wakimoto. "Integrable highest weight modules over affine superalgebras and Appell's function". *Comm. Math. Phys.* 215.3 (2001).

[KW15]   V. G. Kac and M. Wakimoto. *Representations of affine superalgebras and mock theta functions III.* arXiv:1505.01047. 2015.

[LZ14]   S.-C. Lau and J. Zhou. *Modularity of Open Gromov-Witten Potentials of Elliptic Orbifolds.* arXiv:1412.1499. 2014.

[Maa64]  H. Maass. *Lectures on modular functions of one complex variable.* Notes by Sunder Lal. Tata Institute of Fundamental Research Lectures on Mathematics, No. 29. Bombay: Tata Institute of Fundamental Research, 1964.

[Obe12]  G. Oberdieck. *A Serre derivative for even weight Jacobi Forms.* arXiv:1209.5628. 2012.

[Pit09]  A. Pitale. "Jacobi Maaß forms". *Abh. Math. Semin. Univ. Hambg.* 79.1 (2009).

[Ram00]  S. Ramanujan. *Collected papers of Srinivasa Ramanujan.* Edited by G. H. Hardy, P. V. Seshu Aiyar and B. M. Wilson, Third printing of the 1927 original, With a new preface and commentary by Bruce C. Berndt. AMS Chelsea Publishing, Providence, RI, 2000.

[Rau12]  M. Raum. "Dual weights in the theory of harmonic Siegel modular forms". PhD thesis. University of Bonn, 2012.

[Rol15]  L. Rolen. *A new construction of Eisenstein's completion of the Weierstrass zeta function.* arXiv:1504.03787. 2015.

[Sko08]  N.-P. Skoruppa. "Jacobi forms of critical weight and Weil representations". *Modular forms on Schiermonnikoog.* Cambridge Univ. Press, Cambridge, 2008.

[Sko90]  N.-P. Skoruppa. "Developments in the theory of Jacobi forms". *Automorphic functions and their applications (Khabarovsk, 1988).* Khabarovsk: Acad. Sci. USSR Inst. Appl. Math., 1990.

[Vig77]  M.-F. Vignéras. "Séries thêta des formes quadratiques indéfinies". *Séminaire Delange-Pisot-Poitou, 17e année (1975/76), Théorie des nombres: Fasc. 1, Exp. No. 20.* Secrétariat Math., Paris, 1977.

[Zag09]  D. Zagier. "Ramanujan's mock theta functions and their applications (after Zwegers and Ono-Bringmann)". *Astérisque* 326 (2009). Séminaire Bourbaki. Vol. 2007/2008.

[Zag91]  D. Zagier. "Periods of modular forms and Jacobi theta functions". *Invent. Math.* 104.3 (1991).

[Zag94]  D. Zagier. "Modular forms and differential operators". *Proc. Indian Acad. Sci. Math. Sci.* 104.1 (1994). K. G. Ramanathan memorial issue.

[Zwe02]  S. Zwegers. "Mock Theta Functions". PhD thesis. Universiteit Utrecht, 2002.

Max Planck Institute for Mathematics, Vivatsgasse 7, D-53111, Bonn, Germany

E-mail: martin@raum-brothers.eu

Homepage: http://raum-brothers.eu/martin